\documentclass[authoryear,a4paper, 11pt]{elsarticle}

\usepackage{graphicx}
\usepackage{amssymb}

\usepackage[body={8.5in, 11in},left=1 in,right=1 in,top=1 in,bottom=1 in]{geometry}
\usepackage{epsfig,epstopdf}
\usepackage{afterpage, morefloats}
\usepackage{amsmath, booktabs, graphicx, setspace, url, epstopdf}
\usepackage[hang]{caption}
\usepackage{amsthm,amssymb}
\usepackage{amsmath}
\usepackage{mathabx}
\usepackage{algpseudocode}
\usepackage{algorithm,  multirow} 
\usepackage{threeparttable}
\usepackage{natbib}
\usepackage{color}
\usepackage{centernot}
\usepackage{ragged2e} 
\allowdisplaybreaks[1]
\usepackage[titletoc,title]{appendix}
\usepackage{hyperref}
\usepackage{url}
\usepackage{pdflscape}
\usepackage{adjustbox}
\usepackage{longtable}
\usepackage{tabularx}
\usepackage{changepage}
\usepackage{rotating}
\usepackage{verbatim}
\usepackage{lscape}
\usepackage{natbib}
\usepackage{multirow}
\usepackage{hyperref}

\usepackage{lineno}

\journal{}

\begin{document}

\renewcommand{\baselinestretch}{1.5}
\begin{frontmatter}


\title{\textbf{Train timetabling with rolling stock assignment, short-turning and skip-stop strategy for a bidirectional metro line}}



\author[iitd]{Chanchal Kumar Salode}
\ead{smz218543@dms.iitd.ac.in}

\author[iitd]{Prasanna Ramamoorthy\corref{1}}
\ead{prasanna.r@dms.iitd.ac.in}

\address[iitd]{Department of Management Studies, Indian Institute of Technology, Delhi, New Delhi-110016, India}

\cortext[1]{Corresponding author: Prasanna Ramamoorthy, prasanna.r@dms.iitd.ac.in}

\begin{abstract}
Metro train operations is becoming more challenging due to overcrowding and unpredictable irregular passenger demand. To avoid passenger dissatisfaction, metro operators employ various operational strategies to increase the number of train services using limited number of trains. This paper integrates metro timetabling with several operational strategies to improve passenger services with limited number of trains. We propose three optimization models for timetable planning during both peak and off-peak hours: The first model aims to minimize operational costs, the second aims to minimize passenger waiting time, and the third is a multi-objective optimization model that considers both objectives simultaneously. These models integrate operational strategies such as rolling-stock assignment, short-turning, and skip-stop strategies to increase the number of services with limited trains on a bidirectional metro line. The paper also provides detailed calculations for train services, running times, and station dwell times. The proposed models are then implemented on a simplified Santiago metro
line 1.
\end{abstract}

\end{frontmatter}

\setstretch{1.10}

\section{Introduction}\label{sec_intro} 

Rapid transit metro systems play an essential role in city public transportation networks. To satisfy the demand between origin-destination (O-D) pairs, the metro system offers high capacity transportation mode with effective responsiveness and low energy consumption. With increasing urbanization, metro systems are burdened with additional increase in passenger demand that is in excess of the design capacity. Such excess demand, especially during peak hours can cause extreme congestion in the metro stations, leading to passenger dissatisfaction. Under such conditions, one option before the metro operator is to increase the number of trains in service. However, such an approach might not be always possible as it might require investment in terms of additional trains. An alternative approach would be to incorporate some operational strategies in metro line operations to cater to the additional demand with existing number of trains. 

Researchers have typically studied three operational strategies: rolling stock assignment, short-turning, and skip-stop. Rolling stock assignment involves efficiently allocating trains to services, considering service frequency, demand patterns, and cost optimization for enhanced operational performance. Short-turning allows trains to turn back at an intermediate turnaround station before completing their full route, optimizing service efficiency while adapting to varying passenger demand. Beijing's metro implemented a short-turning strategy for 12 metro lines during peak hours amidst the COVID-19 outbreak. This measure aimed to keep train capacity utilization below 50$\%$, thereby lowering the risk of infection (Beijing metro 2020)$\footnote{\url{https://www.sohu.com/a/402944218_380149 }}$. Skip-stop strategy involves specific trains intentionally skipping certain stations to optimize service efficiency and minimize passengers' travel time, achieved by maintaining speed and saving both deceleration and acceleration time along with dwelling time (when the train is stationary at a station). Some lines of the New York City subway adopted skip-stop strategy during peak hours (New York City Subway 2022)$\footnote{\url{https://new.mta.info/}}$. The proper execution of above strategies improves the movement of rolling stock and offers additional train services with existing number of trains. This paper studies a bidirectional metro line with rolling stock assignment, short-turning and skip-stop strategy that increases number of train services in predetermined high-demand operation zones. \\
The main contributions of this work are summarised below. First, the paper outlines a methodology to determine the locations of depots and turnaround facilities on the intermediate station of a metro line. The detailed calculations include various parameters such as potential train services, running time, and dwelling time. Conducting demand analysis enables the metro line to be appropriately divided into multiple operation zones (from one turnaround station to another) to improve the movement of rolling stock and offer additional train services in the high-demand operation zone.

Second, the study conducts a comprehensive analysis of time-varying passenger demand, considering various demand variables for each station pair. During peak hours, some passengers might experience waiting for the second or third train to reach their destination because of the restricted train capacity and predefined operation zones. To represent multiple O-D pairs, along with constrained train capacity and presumed multiple operation zones, the research introduces various demand variables and establishes relationships between them. 

Third, this study is the first to integrate short-turning strategy and rolling stock assignment into metro timetable planning for off-peak hours. For peak hours, in addition to the above strategies, we also include the skip-stop strategy. The incorporation of short-turning and skip-stop strategies enhances rolling stock movement, providing additional train services in the high-demand operation zone. Three optimization models (each for off-peak and peak hours) are formulated: the first, a mixed integer linear programming (MILP) problem minimizing operational costs without considering passenger waiting time; the second, a mixed integer non-linear optimization (MINLP) problem maximizing  service quality by minimizing passengers' waiting time and the third, a multi-objective MINLP
problem that simultaneously minimizing operational costs and passengers' waiting time. The MINLP problem is later transformed into the MILP using linearization techniques. The multi-objective optimization problem is solved using epsilon constraint method, and non-dominating solutions are obtained with respect to both objectives. These models are successfully implemented on the simplified Santiago metro line 1 for validation.

The subsequent sections are structured as follows. Section 2 provides a comprehensive literature review on both regular and irregular timetables, including operational strategies. Section 3 provides a detailed description of the problem and parameter calculations. Section 4 presents model constraints and objective functions. Section 5 introduces mathematical models. Section 6 outlines the implementation of the proposed models on a simplified Santiago metro line 1.  Finally, Section 7 presents the conclusions and directions of future research.
\section{Literature Review}
We present a brief review of the literature from four perspectives: 1) Regular and irregular train timetables; 2) Train timetables with short-turning strategy; 3) Train timetables with rolling stock assignment; 4) Train timetables with skip-stop strategy.
\subsection{Regular and irregular train timetable}
The train timetable problem has attracted considerable attention in recent decades due to its importance in developing metro lines. The research on this problem can be divided into two categories based on metro timetable characteristics: (i) a regular timetable, which allows passengers to easily remember the exact departure times at stations, and (ii) an irregular timetable, which can better accommodate passenger demand without requiring to know the exact departure time while reducing operational costs. However, in this case the time intervals between train arrivals are irregular.

Urban transit metro systems typically employ regular timetables with fixed headways for both peak and off-peak hours. For instance, during peak hours, a train arrives at a station every three minutes, whereas in off-peak hours, the interval increases to every eight minutes. In context of regular train timetable, \cite{yang2015energy}  proposed a scheduling strategy to optimize the use of energy generated by train deceleration. This strategy coordinates the arrival and departure times of trains within the same electricity supply interval, resulting in more efficient energy utilisation. The authors introduced an integer programming model incorporating actual speed profiles and dwell time control to decrease train energy consumption in real-world scenarios.
\cite{canca2017design} proposed a methodology for designing energy-efficient timetables in rapid transit networks. The suggested approach highlights the significance of taking into account all services throughout a planning period and integrating energy consumption into a broader objective function for comprehensive optimization. 

Given the frequent service provided by rapid transit metro systems, passengers often arrive randomly without planning their exact arrival times at the destination (\cite{cepeda2006frequency}). Hence, for dynamically arriving passengers, an irregular timetable is preferable over a regular one. The latter may lead to longer waiting times for passengers and higher operational costs (\cite{barrena2014exact}; \cite{yalccinkaya2009modelling}). 

\cite{cury1980methodology} proposed a methodology which automatically
generates the optimal irregular timetable of metro lines based on train and passenger movements. Similar to model proposed by \cite{cury1980methodology}, \cite{assis2004generation} introduced a predictive control formulation based on linear programming for optimal train schedules. \cite{niu2013optimizing} proposed a binary integer programming
model incorporating passenger loading and departure events to provide a theoretical description of the problem and formulated a non-linear optimization model based on
time-dependent O-D trip records. The authors proposed a local improvement and genetic algorithm for single
and multi-station cases. \cite{canca2014design} presented a MINLP model considering the dynamic behaviour of demand that synchronizes train arrival and departure times. \cite{barrena2014exact} introduced three mathematical formulations to minimize passengers' average waiting time in a dynamic demand scenario. The authors presented a branch-and-cut algorithm applicable to all three formulations. \cite{sun2014demand} addressed the challenges of traditional peak/off-peak schedules in metro services and proposed three optimization models for demand-sensitive timetable design. The authors evaluated models for optimizing metro timetables on a Singapore line, emphasizing the effectiveness of a capacitated, demand-sensitive timetable with dynamic headways. Considering a bi-directional urban metro line, \cite{yin2017dynamic} proposed an integrated approach to address train scheduling problems with the aim of minimizing both operation costs and passengers' waiting time. The authors developed a Lagrangian-based heuristics to find near-optimal solutions in a short computational time. \cite{dong2020integrated} presented
MINLP formulation to optimize train stop planning and timetabling for commuter railways under
time-dependent passenger demand. \cite{yin2021timetable} considered dynamic passenger demand and proposed a mathematical formulation for the coordinated train schedule of a metro network to minimize station crowding during peak hours. To efficiently address large-scale instances, the authors developed a decomposition-based adaptive large neighborhood search (ALNS) algorithm.

\subsection{Train timetabling with short-turning strategy}
In literature, researchers have integrated train timetables with short-turning strategies in various ways, relying on different assumptions and objectives. In the context of bus operations, \cite{delle1998service} presented an optimization framework that considered short-turning strategies for both full-length and short-turn lines, incorporating variable vehicle sizes. \cite{tirachini2011optimal} used demand information to develop a short-turning model within a single bus line for an increasing number of services in the high-demand segment of the line. The authors determined analytical expressions for the optimal value of design variables, capacity of vehicles, and position of short-turn limit stations. \cite{gkiotsalitis2019cost} developed a rule-based approach for short-turning and interlining lines aimed at minimizing operational costs and reducing passenger waiting time.

In the context of metro systems, there is a scarcity of literature on short-turning strategies. Unlike bus line transportation systems, metro lines are independent and have a greater passenger demand and service frequency, which necessitates the study of short-turning in its operations.

\cite{canca2016short} analyzed a short-turning strategy for high-demand operation segment of a line to increase the number of services. The authors formulated a MILP optimization model to determine service offset and turnaround points. \cite{zhang2018short} considered timetabling with short-turning strategy for rapid transit line with multiple depots to improve transport capacity. The authors developed a MINLP optimization model to optimize short-turning and full-length train services with predefined headways. The MINLP formulations are then transformed into MILP by linearization properties and solved using commercial solver. 

In the event of a total obstruction in the rail line when a train can not use the entire track, \cite{ghaemi2018macroscopic} applied short-turning strategy on both sides of disruption. The authors formulated an MILP optimization model to identify the optimal stations and short-turning times. \cite{li2019demand} formulated a MINLP to jointly optimize minimum passenger waiting time, operator costs, and train load balancing. Through real-world examples, the authors presented results from a combination of heterogeneous headways (HH) and short-turning (ST) strategies with those of HH alone, ST alone, and a regular timetable for comprehensive comparison. \cite{blanco2020optimization} presented an MILP model to optimize both operational costs and quality of services by considering varying demands over time, interchange stations, short-turns, and technical train features. The authors
developed a novel matheuristic to solve the problem efficiently for real instances.
\cite{yuan2022integrated} proposed an integrated optimization of a timetabling problem with short-turning strategy and rolling stock assignment to maximize the number of services in the high-demand operation zone. The authors created a hybrid algorithm to achieve high-quality solutions and conducted two case studies, one on a simplified metro line and another on Beijing metro line 6, to confirm the effectiveness of the hybrid algorithm.

\subsection{Train timetabling with rolling stock assignment}
Many researchers have studied the integration of timetabling problems with rolling stock assignments. 
\cite{cadarso2012integration} proposed a model that combines train timetables with rolling stock assignments to satisfy
passenger demand more effectively by adjusting the rolling assignment flow. \cite{hassannayebi2016train} introduced a path-indexed non-linear formulation to minimize the average passenger waiting time in an urban transportation system. The authors applied a Lagrangian relaxation method to relax vehicle turnaround constraints. Consequently, the problem decomposed into multiple sub-problems for each path, leading to more efficient solution of large-scale problems compared to commercial solvers classical method. \cite{yue2017integrated} presented a bi-level programming model using train paths and rolling stock indicators as decision
variables. Upper-level problems minimize passenger waiting time and operation costs, whereas lower-level problems
schedule trains by minimizing the number of infeasible train paths.

\cite{canca2018integrated} introduced a MILP formulation to design rolling stock circulation plans while also determining depot locations and quantities. They proposed a three-phase sequential solution approach. The first phase identifies the minimum number of vehicles for the week, the second identifies daily routes for each line, and the third phase determines weekly circulations and depot locations using a genetic algorithm. \cite{wang2018passenger} integrated train timetable and rolling stock assignment in consideration of time-varying demand of passenger for an urban rail transit line. The authors developed three solution approaches for solving the multi-objective MINLP problem to obtain a train schedule and a rolling stock circulation plan. \cite{liu2020collaborative} introduced a collaborative optimization problem, formulated as a MINLP model, to balance train utilization, strategies for controlling passenger flow, and the number of passengers waiting at platforms. The authors applied a Lagrangian relaxation approach to decompose the problem into two smaller sub-problems, thus reducing the computational burden of the original problem for more efficient solutions. \cite{mo2021exact} developed MINLP model to optimize service quality and operation costs by integrating train timetables and rolling stock assignment. 
The authors decomposed the original model into several sub-problems and solved them using a dynamic programming algorithm. \cite{wang2021real} studied integrated train
rescheduling and rolling stock assignment problems for the complete blockage situations
in a metro line. The authors presented a multi-objective MILP formulation to minimize timetable deviation
and headway deviations of train services.
\subsection{Train timetabling with Skip-stop strategy}
The train timetabling with skip-stop strategy has been studied by many researchers in the literature. \cite{freyss2013continuous} focused on a skip-stop strategy for rail transit lines utilizing a single one-way track. The author introduced a continuous approximation approach to address the problem. In the context of urban rail transit systems, \cite{wang2014efficient} proposed a bi-level optimization approach for solving the MINLP train scheduling problem with a skip-stop strategy. The model aims to simultaneously reduce the total travel time for passengers and minimize the energy consumption of trains. 
For Shanghai-Hangzhou high-speed rail corridor in China, \cite{niu2015train} proposed a non-linear mixed-integer programming model to minimize passenger waiting time using predetermined train skip-stop patterns. \cite{jiang2019q} formulated an integer linear programming problem, representing a collaborative optimization of train scheduling and skip-stop patterns. \cite{zhu2020integrated} presented flexible skip-stop patterns and short-turning techniques for train scheduling problems, considering dynamic passenger behaviour. \cite{chen2022real} investigated the skip-stop strategy to handle real-time disturbances and introduced a non-linear programming model to minimize overall train deviations.

Table 1 compares key features of related literature with our work. Based on the analysis of relevant literature and the data shown in Table 1, it is evident that our paper differs significantly from the literature in two aspects. $(1)$ To the best of our knowledge, this is the first paper that models train timetabling with rolling stock assignment, short-turning and skip-stop strategy. $(2)$  This paper introduces three separate mathematical models for peak and off-peak hours. The first model focuses solely on minimizing operational costs, the second model aims to maximize service quality, while third integrates both as a multi-objective optimization model.
\begin{table}[H]
\caption{Comparison with related literature}
\resizebox{\columnwidth}{!}{%
\begin{tabular}{llllllll}
\hline
Topic & \begin{tabular}[c]{@{}l@{}}Publications\\ (Author (year))\end{tabular} & Research Problem & Model Type & Objective(s) & \begin{tabular}[c]{@{}l@{}}Solution\\ Approach\end{tabular} & \begin{tabular}[c]{@{}l@{}}Network\\ Infrastructure\end{tabular} & \begin{tabular}[c]{@{}l@{}}Train\\ Capacity\end{tabular} \\ \hline
\multirow{7}{*}{Timetable} & Niu \& Zhou (2013) & TT & Non-linear & PWT & LI \& GA & Bi-directional & Yes \\
 & Barrena et al. (2014a) & TT & Non-linear & \begin{tabular}[c]{@{}l@{}}PWT \\ (Three models)\end{tabular} & B\&C & Uni-directional & Yes \\
 & Barrena et al. (2014b) & TT & Non-linear & \begin{tabular}[c]{@{}l@{}}PWT \\ (Two models)\end{tabular} & ALNS & Uni-directional & No \\
 & Canca (2014) & TT & Non-linear & PWT & ALNS & Uni-directional & Yes \\
 & Sun et al. (2014) & TT & Linear & \begin{tabular}[c]{@{}l@{}}PWT \\ (Three models)\end{tabular} & Solver & Uni-directional & Yes \\
 & Yin (2017) & TT & Linear & PWT-EC & LR & Bi-directional & Yes \\
 & Dong (2020) & TT & Non-linear & PWT-DT-RT & ALNS & Bi-directional & Yes \\ \hline
\multirow{18}{*}{Integrated Models} & Canca (2016) & TT-STS & Non-linear & PWT-OC & Solver & Uni-directional & Yes \\
 & Zhang (2018) & TT-STS & Linear & TLV-TSO & Solver & Bi-directional & Yes \\
 & Ghaemi (2018) & TT-STS & Linear & TD & Solver & Network & Yes \\
 & Li (2019) & TT-STS & Non-linear & PWT-OC-TLF & GA & Uni-directional & No \\
 & Blanco (2020) & LP-TT-STS & Linear & PWT-OC & CH & Uni-directional & Yes \\
 & Cardarso (2012) & TT-RSA & Linear & PWT-OC & Solver & Uni-directional & Yes \\
 & Hassannayebi (2016) & TT-RSA & Non-linear & PWT & LR & Bi-directional & Yes \\
 & Yue (2017) & TT-RSA & Linear & WPN & SA & Bi-directional & Yes \\
 & Canca \& Barrena (2018) & TT-RSA & Linear & PWT & GA & Uni-directional & Yes \\
 & Wang (2018) & TT-RSA & Non-linear & \begin{tabular}[c]{@{}l@{}}PWT-TLF-TSO\\ (Multi-objective)\end{tabular} & Iterative-Solver & Bi-directional & Yes \\
 & Liu (2020) & TT-RSA & Linear & WPN-TSO & LR & Bi-directional & Yes \\
 & Mo (2021) & TT-RSA & Non-linear & PWT-OC & DP & Bi-directional & Yes \\
 & Wang (2014) & TT-SSS & Non-linear & RT-EC & Bi-level-Solver & Uni-directional & Yes \\
 & Niu (2015) & TT-SSS & Non-linear & PWT & DP \& GA & Bi-directional & Yes \\
 & Zhu (2020) & TT-SSS & Linear & PWT-GTT & AFOA & Bi-directional & No \\
 & Chen (2022) & TT-SSS & Non-linear & PWT & \begin{tabular}[c]{@{}l@{}}Quadratic Programming\\ based MPC method\end{tabular} & Bi-directional & Yes \\
 & Yuan (2022) & TT-RSA-STS & Non-linear & PWT & GA-Solver & Bi-directional & Yes \\
 & \textbf{This Paper} & \textbf{TT-STS-RSA-SSS} & \textbf{Linear, Non-linear} & \textbf{\begin{tabular}[c]{@{}l@{}}PWT-OC\\ (Three Models\\ (each for peak and \\ off-peak hours))\\ (Multi-objective)\\ (Multiple depots)\end{tabular}} & \textbf{\begin{tabular}[c]{@{}l@{}}Solver-Epsilon constraint\\ method\end{tabular}} & \textbf{Bi-directional} & \textbf{Yes} \\ \hline
\multicolumn{8}{l}{\begin{tabular}[c]{@{}l@{}}Operation strategies: train timetable (TT); line planning (LP); rolling stock assignment (RSA); short-turning strategy (STS); skip-stop strategy (SSS); Objective: passengers' \\ waiting time (PWT); generalized travel time (GTT); operation cost (OC); train loading fluctuation (TLF); energy consumption (EC); delay time (DT); running time \\(RT); train service  connection (TSO); train loading variation (TLV); train delay (TD); waiting passenger number (WPN); Solution approach: general-purpose optimization \\ solver (Solver); genetic algorithm (GA); Lagrangian  relaxation (LR); dynamic programming (DP);  adaptive large neighborhood search algorithm (ALNS); branch and cut \\ (B\&C); local improvement (LI); adaptive fix and optimize algorithm (AFOA); model predictive control (MPC); constructive heuristics (CH); simulated annealing (SA)\end{tabular}}
\end{tabular}
}
\end{table}
\section{Problem Description}\label{sec_prob_descrip}
For a bidirectional metro line with irregular passenger demand, we incorporate the train timetable, rolling stock assignment, short-turning, and skip-stop strategy within four predefined operation zones using the fixed number of trains. 
The line comprises four depots, each located at turnaround stations, and includes 2$N$ stations, as illustrated in Fig. 1. The direction of running from terminal station 1 to station $N$ is termed upstream, whereas the direction from terminal station $N$+1 to station 2$N$ is termed downstream. Two station types are used: a station with a back-turning facility is referred to as a \textit{turnaround} station, and otherwise, it is a \textit{standard} station. In Fig. 1, standard stations are represented by hollow circles, while turnaround stations are denoted by solid circles. Line terminates at turnaround stations (1, $N$, $N$+1, 2$N$) in both directions, serving as the terminals of the metro lines.
\begin{figure}[h!]
\begin{center}
\includegraphics[width=14cm]{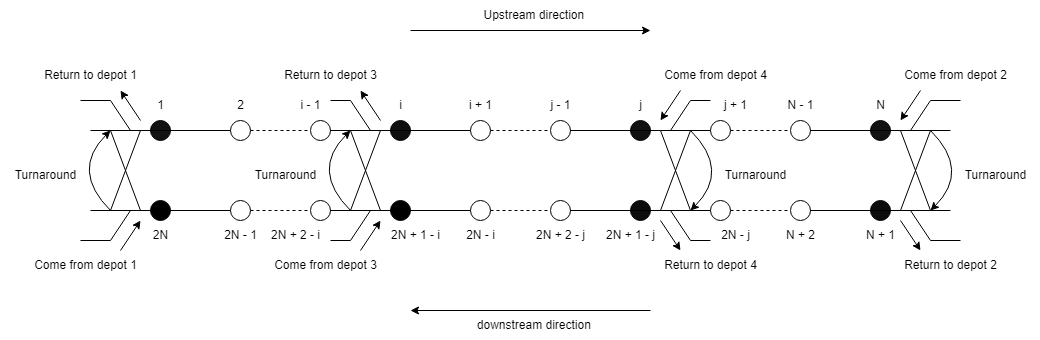}
\end{center}
\caption{Layout of a metro line}
\end{figure} 

We first determine potential services for both upstream and downstream directions through demand and train capacity calculations. The demand matrix indicates the number of passengers traveling from station $i$ to station $j$. By summing the passengers separately for the upstream and downstream directions and dividing by the train carrying capacity, we obtain an upper bound on the services required for the given time interval. Subsequently, we define the upstream and downstream potential service sets $K^{(up)}$ and $K^{(dn)}$, respectively. Service can only begin from one turnaround station, and terminate at another; the segment from one turnaround station to another is known as \textit{operation zone} of the service.

Train service operating between one terminal station (i.e., 1) and another terminal station (i.e., $N$) is termed a \textit{full-length service}. This is illustrated by the services provided by trains in Fig. 2. Conversely, if a train reverse its direction at an intermediate turnaround station, the service is referred to as a \textit{short-turning service}. We illustrate the advantages of a short-turning service in the following example. 

Figures 2, 3, and 4 illustrate scenarios both with and without the short-turning strategy, employing an equal number of trains. To prevent congestion and reduce passenger wait times, the maximum number of possible services in the upstream direction is capped at 8. Assuming increased demand between station 2 and station 4, the full-length strategy in Fig. 2 offers only 6 upstream services. In contrast, train 1 uses a short-turning strategy in Fig. 3, resulting in 8 upstream services. If both trains 2 and 3 adopt a short-turning strategy (see Fig. 4), 9 upstream services can be offered between stations 2 and 4. Thus, when passenger demand varies across segments, a strategically implementing short-turning strategy offers more services on segments with higher passenger demand compared to employing only full-length services, thereby potentially reducing passenger congestion and total waiting time.
\begin{figure}[H]
    \centering    \includegraphics[width=11cm]{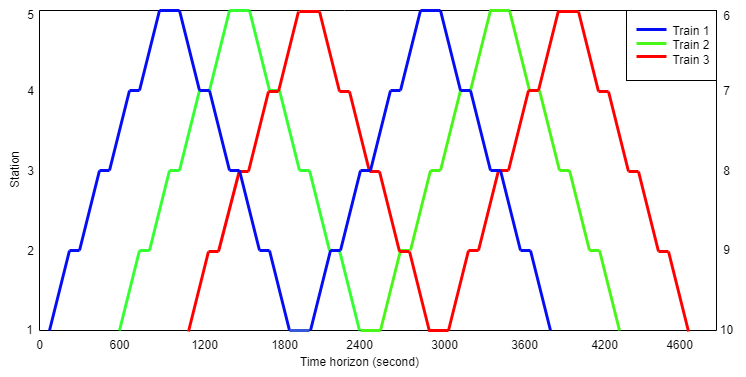}
  \caption{Without short-turning strategy} 
  \end{figure}
  \begin{figure}[H]
    \centering    \includegraphics[width=11cm]{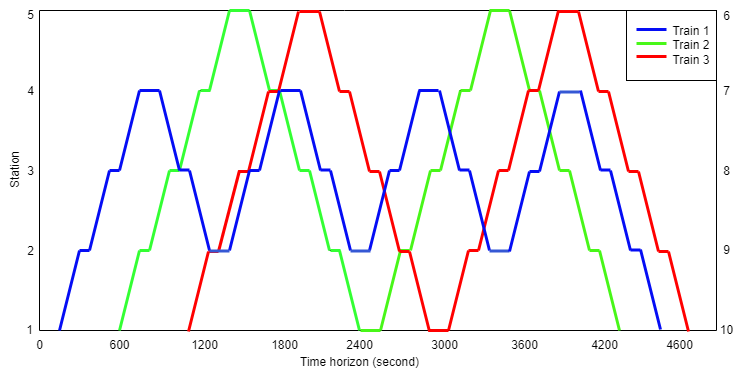}
    \caption{With short-turning strategy (Train 1)} 
  \end{figure}
   \begin{figure}[H]
    \centering    \includegraphics[width=11cm]{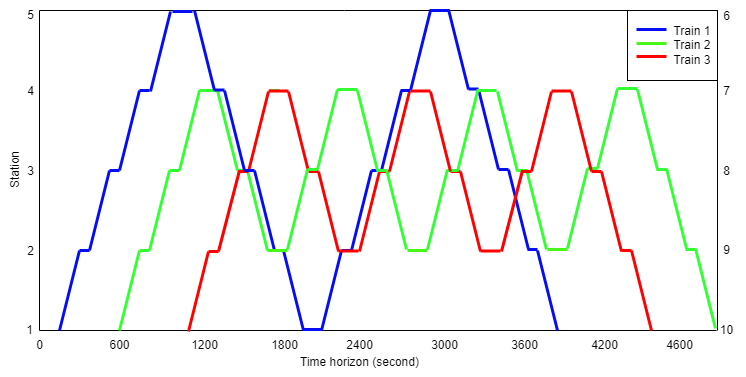}
    \caption{With short-turning strategy (Train 2 and Train 3)} 
  \end{figure}
  
In addition to the circulation of trains and the short-turning strategy, the effective implementation of a skip-stop strategy significantly minimize congestion and passengers' total waiting time. Fig. 5 illustrates an example to explain and compare the skip-stop strategy with strategies in which the train stops at all stations (we call it standard stop strategy).

Assuming there are four stations with 200, 500, and 200 passengers at stations 1, 2, and 3, respectively, all headed to station 4. The train can accommodate up to 600 passengers. Generally, the rate of new passenger arrivals is much lower than the number of passengers already waiting at the stations.

In Fig. 5(a), the standard stopping strategy is employed. Service $u_{1}$ originates from station 1 with 200 passengers and proceeds to stations 2 and 3, boarding 400 passengers from station 2. However, due to capacity restrictions, 100 and 200 passengers remain waiting at station 2 and station 3, respectively. Service $u_{2}$ departs from stations 2 and 3 with the remaining 100 and 200 passengers, respectively, completing its total travel time in $T_{std}$ = 13. The calculation of total passenger waiting time in Fig. 5(a) is as follows:
\begin{align*}
&Wt_{std} = 200 \times 2 + 500 \times 5 + 100 \times 3 + 200 \times 11 = 5400.
\end{align*}

The skip-stop strategy is employed in Fig. 5(b). Service $u_{1}$ skips station 1 and departs from stations 2 and 3 with 500 and 100 passengers, respectively, leaving 100 passengers at station 3 due to capacity restrictions. Service $u_{2}$ departs from station 1, skips station 2, and boards 100 passengers from station 3. The travel time finishes at $T_{skp}$ = 12. The calculation of passengers' waiting time using the skip-stop strategy is as follows:
\begin{align*}
&Wt_{skp} = 200 \times 5 + 500 \times 4 + 200 \times 7 + 100 \times 3 = 4700.
\end{align*}
By implementing the skip-stop strategy, both the finishing time and passenger waiting time are effectively reduced by:
\begin{align*}
&\frac{T_{std}-T_{skp}}{T_{std}} = \frac{13-12}{13} = 7.69\% \\
&\frac{Wt_{std}-Wt_{skp}}{Wt_{std}} = \frac{5400-4700}{5400} = 12.96\%
\end{align*}
\begin{figure}[H]
    \centering    \includegraphics[width=9cm]{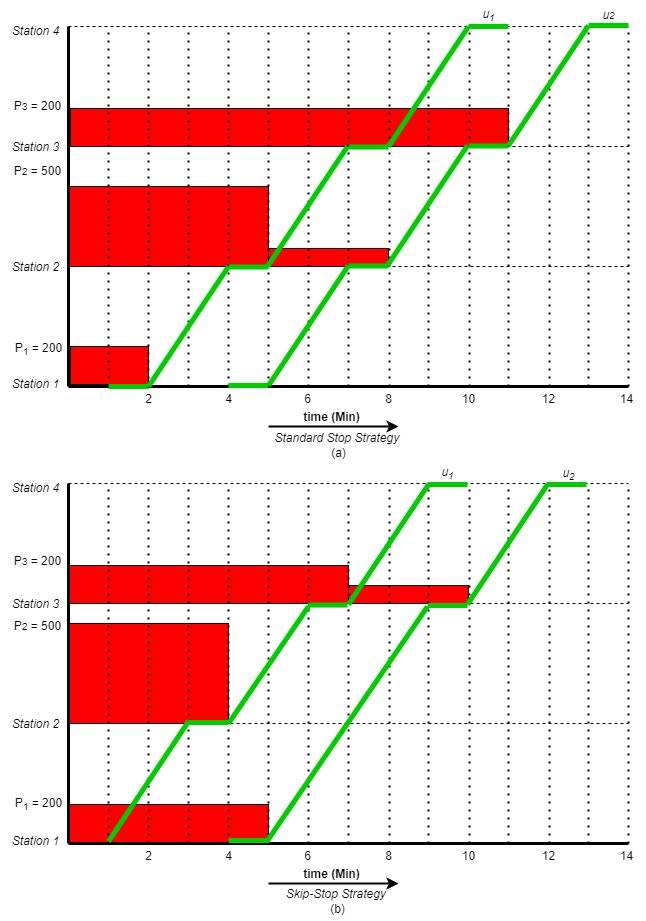}
  \caption{Without short-turning strategy} 
  \end{figure}
As per the example mentioned above, it is clear that this approach significantly reduces travelling time and passenger waiting time during peak hours. Moreover, the skip-stop strategy enables better resource utilization, minimizes congestion and optimizes service efficiency. This results in an improved overall passenger experience, making it an effective approach for managing peak-hour demands in metro systems.

In summary, this paper  addresses four interconnected decision problems. The first involves selecting service operation zones from predefined operation zones for each service. The second decision is to create a train schedule by determining the arrival and departure times for each service. The third decision involves allocating rolling stock and a limited number of trains to services, ensuring that each service is carried out by a single train. The fourth decision integrates the train timetable, assignment of rolling stock, and short-turning strategy for off-peak hours, along with the skip-stop strategy for peak hours. These four decision problems are simultaneously addressed in the optimization models described in Section 5.
\noindent The mathematical models developed in this study are based on the following assumptions.
\\
\textbf{Assumption 1:} There are four depots; two are located at the terminals, and two at intermediate stations. The identification of locations for intermediate depots is discussed in subsection 3.1.\\ 
\textbf{Assumption 2:} Running time between station pairs and dwelling time at each station are predetermined, as explained in subsections 3.3 and 3.4, respectively.\\
\textbf{Assumption 3:} During off-peak hours, the train timetable, rolling stock assignment and short-turning strategy are integrated, while peak hours involve the integration of train timetable, short-turning, rolling stock assignment, and skip-stop strategy. \\
\textbf{Assumption 4:} No train overtaking, connections, or meetings are permitted.\\
\textbf{Assumption 5:} Passengers are restricted from transferring between services; they may only board a service that directly reaches their destination station.\\
\textbf{Assumption 6:} 
Passengers board according to the ``first-in, first-out'' (FIFO) approach.\\
\\
\textbf{Assumption 1} 
specifies location and number of depots for a line, with two depots at each terminal station and the other two at intermediate turnaround stations, identified through the process explained in subsection 3.1. According to \textbf{Assumption 2}, running and dwell times are predetermined through appropriate calculations to operate services (explained in section 3.3 and 3.4, respectively). \textbf{Assumption 3} outlines integration of train timetable, rolling stock assignment and short-turning strategy for off-peak hours, and skip-stop strategy in addition to the above strategies is considered for peak hours. Passengers are informed in advance about which stations the next services will skip during peak hours. Additionally, the train (in-service) can adopt a skip-stop strategy after departing from the end station of the service operation zone.  \textbf{Assumption 4} states that each station is served by a single train at a time, with no meeting, overtaking, or train connections allowed in this study. \textbf{Assumption 5}
assumes passengers can only board the service that travels directly to their destination. The exact boarding order of passengers is determined by \textbf{Assumption 6}.
\subsection{Locating depots on intermediate stations}
To determine the location of depots at intermediate stations, it is essential to identify a continuous segment of stations with the highest passenger flow compared to other segments. This involves selecting segments situated between stations $i-1$ and $j+1$ (as illustrated in Fig. 1) and subsequently subdividing all stations into sub-segments, such as $i-1$ to $i+4$, $i-1$ to $i+5$, ..., $i-1$ to $j-1$, $i-1$ to $j$, $i-1$ to $j+1$, with the assumption that each sub-segment contains at least four stations. The segment with the highest demand is then chosen from all sub-segments, and depots are established at intermediate stations, with turnaround facilities positioned at both the starting and ending stations of the selected segment.
\subsection{Train Service Calculations}
The operational period for an urban transit line is defined as [$t_{\text{start}}$, $t_{\text{end}}$]. This period is divided into multiple time intervals, denoted as [$t_{0}$, $t_{1}$], [$t_{1}$, $t_{2}$], ..., [$t_{s}$, $t_{s+1}$], ..., [$t_{S-1}$, $t_{S}$] where $t_{0}$ = $t_{\text{start}}$ and $t_{S}$ = $t_{\text{end}}$. The number of potential train services for these time intervals is denoted by $K_{1}$, $K_{2}$, ..., $K_{s}$, $K_{s +1}$, ..., $K_{S}$, which is calculated based on the total passenger demand between stations for the interval [$t_{s-1}$, $t_{s}$]:
\begin{align*}
&\Omega_{s} = C_{train}.\kappa . K_{s} &\forall s \in \{1,2,....,S\}
\end{align*}
Where $\Omega_{s}$, $\kappa$, and $C_{\text{train}}$ represent the passenger demand for the time interval [$t_{s-1}$, $t_{s}$], the predefined load factor (loading capacity of train), and the train's carrying capacity, respectively. The predefined load factor is employed to limit overcrowding of onboard passengers and is typically set to a value less than 1. During off-peak hours, a common value is 0.70, indicating that if the train has a capacity of 1000 passengers (sitting and standing), the maximum number of passengers per train is 700. Similarly, during peak hours, it could reach 0.95 or even 1.
\subsection{Running time calculations}
Running time between two stations $i$ and $i+1$ consists of three components: acceleration, deceleration, and pure running time (as illustrated in Fig. 6). To calculate the total time required to cover the distance between two consecutive stations, the following equation can be utilized:
\begin{align*}
Running\_time =&Time_{acc} + Time_{dec} + Time_{pure}&
\end{align*}
Let $V_{acc}$, $V_{dec}$, and $V_{max}$ represent the acceleration, deceleration, and maximum speed of a train, respectively. Where $Time_{acc}$ is the time required for acceleration, calculated as $Time_{acc} = \frac{V_{max}}{V_{acc}}$, $Time_{dec}$ is the time required for deceleration, calculated as $Time_{dec} = \frac{V_{max}}{V_{dec}}$, and $Time_{pure}$ is the time between consecutive stations without considering acceleration and deceleration time, calculated using the equation: $Time_{pure} = \frac{\text{(Remaining distance)}}{V_{max}}$. 
\begin{figure}[h!]
\begin{center}
\includegraphics[width=10cm]{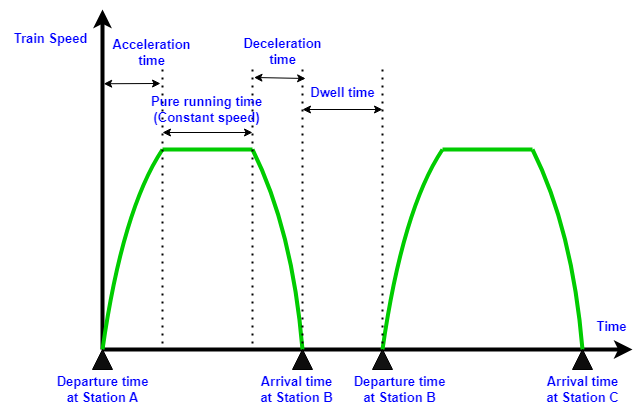}
\end{center}
\caption{Train time-speed representation}
\end{figure}

\subsection{Dwell time calculations}
The time difference between the departure and arrival of a train at station B (refer to Fig. 6) is the dwell time of the station, which directly depends on the number of arrivals and departures (crowdedness) at a station. 

Let $p_{i,j}$ denote the O-D demand matrix representing the number of commuters traveling from station $i$ to $j$. The dimension of matrix  $p_{i,j}$ is $N$ $\times$ $N$, where $S$ and $N$ represents the set and total number of stations, respectively.

The passenger concentration at station $s$ can be calculated at particular time interval as follows:
\begin{align*}
&Crowdedness_{s} = \Bigr[\sum_{i}p_{i,s} + \sum_{j}p_{s,j}\Bigr]& \forall s \in S     
\end{align*}
In the above expression, the first term represents the total number of passengers arriving at station $s$ from all other stations $i$, and the second term indicates the total number of passengers departing from station $s$ to all stations 
$j$. The sum of these two terms defines the crowdedness of station $s$.

The crowdedness of any station determines the dwell time for boarding and alighting of passengers. Hence, based on the crowdedness values at stations, different groups of stations are formed using a predefined threshold, and specific dwell times are allocated to each group. For instance, if the crowdedness of station 1, 2, and 3 surpasses certain thresholds, they are categorized into different groups, and respective dwell times, such as 40, 45, and 35 seconds, are allocated to each group. 


\section{Constraints and Objective functions}
\noindent Table 2 represents all the symbols and parameters used in the mathematical formulations.
\begin{table}[htb]
\caption{Model Parameters}
\resizebox{\columnwidth}{!}{%
\begin{tabular}{ll}
\hline
Parameters & Detailed Description                                                                     \\ \hline
$S^{(up)}$      & Set of stations in the upstream direction, $S^{(up)}$ $=$ \{1 , . . . , N\}                     \\
$S^{(dn)}$       & Set of stations in the downstream direction, $S^{(dn)}$ $=$ \{N+1 , . . . , 2N\}                  \\
$S^{(up)}_{M}$      & Set of turnaround stations in the upstream direction, i.e., $S^{(up)}_{M}$$\subset$ $S^{(up)}$                  \\
$S^{(dn)}_{M}$      & Set of turnaround stations in the downstream direction, i.e., $S^{(dn)}_{M}$$\subset$ $S^{(dn)} $               \\
$K^{(up)}$      & Set of potential train services in upstream direction                                           \\
$K^{(dn)}$        & Set of potential train services in downstream direction   
 \\
$D$        & Set of depots          \\
$i,j,m,n$    & Index of stations                                                                            \\
$k,l$        & Index of train services                                                                           \\
$d$        & Depot index  
 \\
$s_{k}$        & Maximum number of stations that service $k$ can skip
\\
$h_{max}$       & Maximum headway between two subsequent services                                        \\
$h_{min}$        & Minimum headway between two subsequent services 
\\
$Rt_{j}$   & Pure running time from station $j$ $-$ $1$ to station $j$         \\
$R_{a}$   & Extra acceleration time
\\
$R_{d}$   & Extra deceleration time
\\
$r_{j-1,j}$   & Running time from station $j$ $-$ $1$ to station $j$ ($Rt_{j}$ + $R_{a}$ + $R_{d}$)                           \\
$e_{i}$         & Dwelling time at station $i$                                               \\
$\delta^{min}_{m}$ & Minimum back-turning time at station $m$ (turnaround station)                                       \\
RS         & Total number of trains available for operations              \\
$P_{ij}[t]$ & Arrival rate of passengers travelling from station $i$ to destination station $j$ at time $t$ \\
C          & Capacity of a train \\
M          & A large positive number                                               \\ \hline
\end{tabular}
}
\end{table}

\subsection{Decision variables}
Table 3 presents decision variables for selection of operation zones, the design of train timetables, and assignment of rolling stock.
\begin{table}[H]
\caption{Decision variables}
\resizebox{\columnwidth}{!}{%
\begin{tabular}{ll}
\hline
Decision Variables             & Definition                                \\ \hline
$\tau_{k}$               & \begin{tabular}[c]{@{}l@{}}0-1 Binary variable; equal to 1 , if service $k$ is selected from potential services,\\ 0, otherwise\end{tabular} \\
$z_{k,m,n}$               & \begin{tabular}[c]{@{}l@{}}0-1 Binary variable; equal to 1 , if service $k$ start from turnaround station $m$ and \\ end  at turnaround station $n$ respectively, and $m\textless{}n$, 0, otherwise\end{tabular} \\
$y_{k,l,m}$           & \begin{tabular}[c]{@{}l@{}}0-1 Binary variable; equal to 1 , when train finishes service $k$ at turnaround station \\ $m$ and starts service $l$ in opposite direction, 0, otherwise\end{tabular}                 \\
$x_{k,i}$            & 0-1 Binary variable; equal to 1 , if station $i$ served by service $k$, 0, otherwise\\
$\alpha^{(up,dp)}_{k}$               & \begin{tabular}[c]{@{}l@{}}0-1 Binary variable; equal to 1 , if service $k$ in the upstream direction is performed \\ by a train arriving directly from depot $dp$, 0 otherwise \end{tabular} \\
$\beta^{(up,dp)}_{k}$               & \begin{tabular}[c]{@{}l@{}}0-1 Binary variable; equal to 1 , if a train returns to depot $dp$ after completing \\ upstream service $k$ , 0, otherwise\end{tabular} \\
$\alpha^{(dn,dp)}_{k}$               & \begin{tabular}[c]{@{}l@{}}0-1 Binary variable; equal to 1 , if service $k$ in the downstream direction is performed \\ by a train arriving directly from depot $dp$, 0, otherwise\end{tabular} \\
$\beta^{(dn,dp)}_{k}$               & \begin{tabular}[c]{@{}l@{}}0-1 Binary variable; equal to 1 , if a train returns to depot $dp$ after completing \\ downstream service $k$ , 0, otherwise\end{tabular} \\
$h_{l-1,l}$             & Headway between two consecutive services ($l-1$ and service $l$) \\
$a_{k,i}$                  & Arrival time at station $i$ of train service $k$  \\
$d_{k,i}$                  & Departure time at station $i$ of train service $k$\\
$RS_{dp}$ &  The number of trains assigned to depot $dp$ before operations commence \\
\hline
\end{tabular}
}
\end{table}

\subsection{Constraints} According to Assumption 3, the integration of train timetable, rolling stock assignment, and short-turning strategy is taken into account for off-peak hour operations. Thus, there are slight variations in constraints between off-peak and peak hour operational planning. The operational constraints specific to off-peak hours are listed below.
\subsubsection{Off-peak hours}
The mathematical model comprises four distinct categories of constraints: 1) Operation zone selection constraints, 2) Train timetable constraints, 3) Rolling stock assignment constraints, and 4) Passenger demand constraints.\\
\textbf{Operation zone constraints:}\\
\textbf{(i) Operation zone selection constraints:}\\
The identification of depot locations at intermediate stations is determined through the process explained in Section 3.1. For simplicity, we assume that the traffic between station 3 to station $N$-2 or stations $N$+3 to 2$N$-2 is high compared to other segments. As a result, depots 3 and 4 are located at stations 3 (upstream) or 2$N$-2 (downstream) and $N$-2 (upstream) or $N$+3 (downstream), respectively. Consequently, more services are required in this segment, and this can be achieved by employing a short-turning strategy. Thus, the operation zones for upstream services are divided into four segments (refer Fig. 1): station 1 to station $N$, station 1 to station $N$-2, station 3 to station $N$, and station 3 to station $N$-2. Similarly, for downstream services: station $N$+1 to station 2$N$, station $N$+1 to station 2$N$-2, station $N$+3 to station 2$N$, and station $N$+3 to station 2$N$-2.
\begin{align}
&\sum_{m\in S^{(up)}_M/\{N-2\}} \sum_{n\in S^{(up)}_M/\{3\}}z_{k,m,n} = \tau_{k} & \forall k \in K^{(up)} \label{1}\\
&\sum_{m\in S^{(dn)}_M/\{2N-2\}} \sum_{n\in S^{(dn)}_M/\{N+3\}}z_{k,m,n} = \tau_{k} & \forall k \in K^{(dn)} \label{2}
\end{align}
Constraints (\ref{1}) and (\ref{2}) restrict the number of potential operation zones to four sub-segments for each stream. These constraints specify that if a potential service $k$ is chosen, then one operation zone for the upstream direction and one for the downstream direction are assigned to that service.\\
\textbf{(ii) Served station constraints}\\
If a potential service $k$ belonging to the set $K^{(up)}$ is chosen (i.e., $\tau_{k} = 1$), then only service $k$ can serve station $i$ within the set $S^{(up)}$ (i.e., $x_{k,i} = 1$ for $i \in S^{(up)}$).
\begin{align}
&x_{k,i} \leq \tau_{k} & \forall k \in K^{(up)}, i \in S^{(up)}\label{3}
\end{align}
If a potential upstream service $k$ is chosen (i.e., $\tau_{k}=1$) to serve operation zone $z_{k,m,n}$ (i.e., $z_{k,m,n} = 1$), then it cannot serve any station beyond station $m$ (i.e., $\sum_{i<m}x_{k,i} = 0$) or after station $n$ (i.e., $\sum_{i>n}x_{k,i} = 0$). This constraint ensures this limitation.
\begin{align}
&\sum_{i<m} x_{k,i} \leq M.(1 - \sum_{n\in S^{(up)}_{M}}z_{k,m,n}) & \forall k \in K^{(up)}, m\in S^{(up)}_{M}/\{1,N-2\}, i \in S^{(up)}\label{4}\\
&\sum_{i>n} x_{k,i} \leq M.(1 - \sum_{m\in S^{(up)}_{M}/ \{N-2\}}z_{k,m,n}) & \forall k \in K^{(up)}, n\in S^{(up)}_{M}/\{3\}, i \in S^{(up)}\label{5}
\end{align}
Moreover, upstream service $k$ must provide service to all stations within its designated operation zone, as ensured by the following constraints.
\begin{align}
&x_{k,i} \geq z_{k,m,n} & \forall k \in K^{(up)}, i \in S^{(up)},m \in S^{(up)}_{M}/\{N-2\},n \in S^{(up)}_{M}/\{3\}, m\leq i \leq n \label{6}
\end{align}
Constraints (\ref{3})-(\ref{6}) establish the relationship between variables $x_{k,i}$, $z_{k,m,n}$, and $\tau_{k}$ for potential upstream service $k \in K^{(up)}$. Likewise, for potential downstream service $k \in K^{(dn)}$, the following constraints need to be satisfied:
\begin{align}
&x_{k,i} \leq \tau_{k} & \forall k \in K^{(dn)}, i \in S^{(dn)}\label{7}
\end{align}
\begin{align}
&\sum_{i<m} x_{k,i} \leq M.(1 - \sum_{n\in S^{(dn)}_{M}} z_{k,m,n}) & \forall k \in K^{(dn)}, m\in S^{(dn)}_{M}/\{N+1,2N-2\}, i \in S^{(dn)}\label{8}\\
&\sum_{i>n} x_{k,i} \leq M.(1 - \sum_{m\in S^{dn}_{M}/ \{2N-2\}}z_{k,m,n}) & \forall k \in K^{(dn)}, n\in S^{(dn)}_{M}/\{N+3\}, i \in S^{(dn)}\label{9}
\end{align}
\begin{align}
&x_{k,i} \geq z_{k,m,n} & \forall k \in K^{(dn)}, i \in S^{(dn)},m \in S^{(dn)}_{M}/\{2N-2\},n \in S^{(dn)}_{M}/\{N+3\}, m\leq i \leq n \label{10}
\end{align}
Moreover, according to Assumption 4, passengers are only permitted to board a service that caters to both their origin and destination. Therefore, an appropriate selection of the operation zone must be made; otherwise, certain passengers might face prolonged waiting times. To ensure this condition, station $i$ must be serviced by at least one service in every pair of consecutive services in each direction, i.e.,
\begin{align}
&x_{k-1,i} + x_{k,i} \geq 1 & \forall k-1,k \in K^{(up)},i \in S^{(up)}\label{11}\\
&x_{k-1,i} + x_{k,i} \geq 1 & \forall k-1,k \in K^{(dn)},i \in S^{(dn)}\label{12}
\end{align}
\textbf{Train timetable constraints:}\\
\textbf{(i) Arrival and Departure time constraints:}\\
As per Assumption 2, the running and dwell times are predetermined through time-speed and crowdedness calculations, respectively (explained in Section 3.3 and 3.4). 
The following constraints provided below define when potential service $k$ arrives and departs from each station.
\begin{align}
&a_{k,i} - d_{k,i-1} = r_{i-1,i} & \forall k \in K^{(up)}, i-1,i \in S^{(up)} \label{13}\\
&a_{k,i} - d_{k,i-1} = r_{i-1,i} & \forall k \in K^{(dn)}, i-1,i \in S^{(dn)} \label{14}\\
&d_{k,i} - a_{k,i} = e_{i} & \forall k \in K^{(up)}, i \in S^{(up)} \label{15}\\
&d_{k,i} - a_{k,i} = e_{i} & \forall k \in K^{(dn)}, i \in S^{(dn)} \label{16}
\end{align}
To define the time horizon, the departure time from the first station for the first train service in both directions is given ($Ft^{up}_{t}$, for e.g. $Ft^{up}_{t}$ = 25200 (7.00am)), and similarly, the bounds on the departure time of the last train service ($Lt^{up}_{t}$, for e.g. $Lt^{up}_{t}$ = 28800 (8.00am)) at the turnaround station are fixed for all train services. The following
constraints ensure these conditions:
\begin{align}
&d_{1,1} = Ft^{up}_{t}\label{17}\\
&d_{N+1,1} = Ft^{dn}_{t}\label{18}\\
&d_{k,1} \leq Lt^{up}_{t}  & \forall k \in K^{(up)}\label{19} \\
&d_{k,3} \leq Lt^{up}_{t}  & \forall k \in K^{(up)} \label{20} \\
&d_{k,N+1} \leq Lt^{dn}_{t}  & \forall k \in K^{(dn)} \label{21}\\
&d_{k,N+3} \leq Lt^{dn}_{t}  & \forall k \in K^{(dn)}
\label{22}
\end{align}
\textbf{(ii) Headway constraints:}\\
To ensure operational safety and maintain a specified service level, minimum and maximum headway are defined.  If potential service $k$ is selected ($\tau_{k} = 1$), the time interval between $k-1$ and $k$ must fall within the specified minimum and maximum limits. Conversely, if potential service $k$ is not chosen ($\tau_{k}=0$), the departure times of service $k-1$ and $k$ at any station should be the same. The subsequent constraints ensure these maximum and minimum interval limits.
\begin{align}
&h_{min}.\tau_{k} \leq h_{k-1,k} \leq h_{max}.\tau_{k}  & \forall k-1,k \in K^{(up)} \label{23}\\
&h_{min}.\tau_{k}  \leq h_{k-1,k} \leq h_{max}.\tau_{k}  & \forall k-1,k \in K^{(dn)} \label{24}
\end{align}
The headway between service $k-1$ and $k$ at station $i$ are determined by difference between departure times of service $k$ and $k-1$ at station $i$.
\begin{align}
&d_{k,i} = d_{k-1,i} + h_{k-1,k} & \forall k-1,k \in K^{(up)}, i \in S^{(up)} \label{25}\\
&d_{k,i} = d_{k-1,i} + h_{k-1,k} & \forall k-1,k \in K^{(dn)}, i \in S^{(dn)} \label{26}
\end{align} 
\textbf{Turnaround constraints:}\\
After selecting the operation zone for potential services, the turnaround stations are automatically designated as the end stations of that service in both the upstream and downstream directions. These stations have back-turning facilities, enabling trains to reverse their direction from upstream to downstream or vice versa.
To allow a train to transition from providing upstream service $k$ to downstream service $l$, the terminal station of service $k$ ($N$ or $N-2$) and the starting station of downstream service $l$ ($N+1$ or $N+3$) must form a station pair (as illustrated in Fig. 1). The following constraints need to be met for the variable $y_{k,l,m}$:
\begin{align}
&2.y_{k,l,N} \leq z_{k,1,N} + z_{l,N+1,2N} + z_{k,3,N} + z_{l,N+1,2N-2} & \forall k \in K^{(up)},l \in K^{(dn)} \label{27}\\
&2.y_{k,l,N-2} \leq z_{k,1,N-2} + z_{l,N+3,2N} + z_{k,3,N-2} + z_{l,N+3,2N-2}  & \forall k \in K^{(up)},l \in K^{(dn)} \label{28}
\end{align}

Constraints (\ref{27}) and (\ref{28}) ensure that the turnaround variable $y_{k,l,N}$ is set to 1 if the operation zone for upstream service $k$ is selected ($z_{k,1,N} = 1$) and simultaneously the operation zone for downstream service $l$ is chosen ($z_{l,N+1,2N} = 1$). In such a case, the train can turn around at turnaround station $m$, indicated by $y_{k,l,N} = 1$.

If a train turns around at the end station of the operation zone after finishing upstream services $k$ to begin downstream services $l$ (i.e., $y_{k,l,N} = 1$), the departure time for upstream service $k$ at the turnaround station $N$ or $N-2$, and the arrival time for downstream service $l$ at station $N+1$ or $N+3$, must meet the minimum turnaround times.
\begin{align*}
& a_{l,N+1} - d_{k,N} \geq \delta^{min}_{N}&  \text{if }  y_{k,l,N} = 1\\
&a_{l,N+3} - d_{k,N-2} \geq \delta^{min}_{N-2} & \text{if }  y_{k,l,N-2} = 1     
\end{align*}
The following linear constraints are equivalent to the above constraints:
\begin{align}
&a_{l,N+1} - d_{k,N} \geq \delta^{min}_{N} + M.(y_{k,l,N}-1)  & \forall k \in K^{(up)},l \in K^{(dn)} \label{29}\\
&a_{l,N+3} - d_{k,N-2} \geq \delta^{min}_{N-2} + M.(y_{k,l,N-2}-1)  & \forall k \in K^{(up)},l \in K^{(dn)} \label{30}
\end{align}
Similarly, if a train completes downstream service $k$ and then proceeds to upstream service $l$, the following constraints must be satisfied:
\begin{align}
&2.y_{k,l,2N} \leq z_{k,N+1,2N} + z_{l,1,N} + z_{k,N+3,2N} + z_{l,1,N-2}   & \forall k \in K^{(dn)},l \in K^{(up)} \label{31}\\
&2.y_{k,l,2N-2} \leq z_{k,N+1,2N} + z_{l,1,N-2} + z_{k,N+3,2N-2} + z_{l,3,N-2}  & \forall k \in K^{(dn)},l \in K^{(up)} \label{32}\\
&a_{l,1} - d_{k,2N} \geq \delta^{min}_{2N} + M.(y_{k,l,2N}-1)  & \forall k \in K^{(dn)},l \in K^{(up)} \label{33}\\
&a_{l,3} - d_{k,2N-2} \geq \delta^{min}_{2N-2} + M. (y_{k,l,2N-2}-1)  & \forall k \in K^{(dn)},l \in K^{(up)} \label{34}
\end{align}
\textbf{Rolling stock constraints:}\\
The upstream or downstream services are either performed by a new train coming from the depot or by up/downstream services that reverse their direction at the turnaround station. Both possibilities are ensured by the following constraints
\begin{align}
    &\sum_{l \in K^{(dn)}} y_{l,k,2N} + \sum_{l \in K^{(dn)}} y_{l,k,2N-2} + \alpha^{(u,1)}_{k} + \alpha^{(u,3)}_{k} = \tau_{k} & \forall k \in K^{(up)}\label{35}\\
&\sum_{l \in K^{(up)}} y_{l,k,N-2} + \sum_{l \in K^{(up)}} y_{l,k,N} + \alpha^{(d,2)}_{k} + \alpha^{(d,4)}_{k}= \tau_{k} & \forall k \in K^{(dn)}\label{36}
\end{align}
After finishing upstream/downstream services, the train can either return to the depot or turn around to perform downstream/upstream services in the opposite direction. The following constraints ensure these conditions:
\begin{align}
&\sum_{l \in K^{(up)}} y_{k,l,2N} + \sum_{l \in K^{(up)}} y_{k,l,2N-2} + \beta^{(dn,1)}_{k} + \beta^{(dn,3)}_{k} = \tau_{k} & \forall k \in K^{(dn)}\label{37}\\
&\sum_{l \in K^{(dn)}} y_{k,l,N-2} + \sum_{l \in K^{(dn)}} y_{k,l,N}  + \beta^{(up,2)}_{k} + \beta^{(up,4)}_{k} = \tau_{k} & \forall k \in K^{(up)}\label{38}
\end{align}
The subsequent constraints ensure that the overall train allocation to depots adheres to the maximum number of trains available in the network.
\begin{align}
&RS_{1} + RS_{2} + RS_{3} + RS_{4} \leq RS &\label{39}
\end{align}
In this study, we consider four depots, with two situated at terminal stations and the other two at intermediate turnaround stations. Variables $RS_{1}$, $RS_{2}$, $RS_{3}$, and $RS_{4}$ denote the number of trains assigned to depots 1, 2, 3, and 4, respectively, prior to operations. According to the definition of $\alpha^{(up,dp)}_{k}$, the subsequent constraint must be met:
\begin{align}
&RS_{1} \geq \sum_{k \in K^{(up)}} \alpha^{(up,1)}_{k} & \label{40}
\end{align}
This implies that the number of trains originating from depot 1 to conduct upstream services does not exceed the number of trains allocated to depot 1. The same principle applies to depots 2, 3, and 4.
\begin{align}
&RS_{2} \geq \sum_{k \in K^{(dn)}} \alpha^{(dn,2)}_{k} & \label{41}\\
&RS_{3} \geq \sum_{k \in K^{(up)}} \alpha^{(up,3)}_{k} & \label{42}\\
&RS_{4} \geq \sum_{k \in K^{(dn)}} \alpha^{(dn,4)}_{k} & \label{43}
\end{align}
To avoid trains traveling without passengers, a train arriving from Depot 1 must begin its service from Station 1. Likewise, if a train completes an upstream service and returns to Depot 2, the final station for this service is designated as Station $N$. These requirements are outlined as follows:
\begin{align}
&\alpha^{(up,1)}_{k} \leq z_{k,1,N-2} + z_{k,1,N}  & \forall k \in K^{(up)}\label{44}\\
&\beta^{(up,2)}_{k} \leq z_{k,1,N} + z_{k,3,N}  & \forall k \in K^{(up)}\label{45}
\end{align}
Similarly, for downstream service $k$:
\begin{align}
&\alpha^{(dn,2)}_{k} \leq z_{k,N+1,2N-2} + z_{k,N+1,2N}  &\forall k \in K^{(dn)} \label{46}\\
&\beta^{(dn,1)}_{k} \leq z_{k,N+1,2N} + z_{k,N+3,2N}  & \forall k \in K^{(dn)}\label{47}
\end{align}
\textbf{Passenger demand constraints:}\\
In this section, we formulate constraints representing passenger demand on a bidirectional metro line.\\
\textbf{Decision Variables for passenger demand:}\\ 
In this study, we represent the number of passengers traveling from station $i$ to station $j$ at time $t$ by $P_{i,j}[t]$. The variables listed in Table 4, in addition to the operation zone selection variable $x_{k,i}$ and timetable variables $d_{k,i}$, play a role in determining the progression of passenger demand over time.
\begin{table}[h]
\caption{Passenger demand variables}
\resizebox{\columnwidth}{!}{%
\begin{tabular}{ll}
\hline
Decision Variables & Definition                                                                                                                                                                                                    \\ \hline
$w_{k,i,j}$     & Number of passengers waiting at station $i$ with destination $j$ for service $k$                                                                                                                                    \\
$w^{b}_{k,i}$      & Number of passengers who want to board service $k$ at station $i$                                                                                                                                                 \\
$w^{b}_{k,i,j}$     & Number of passengers who want to board service $k$ at station $i$ with destination $j$                                                                                                                              \\
$n^{b}_{k,i}$     & Number of passengers who can board service $k$ at station $i$                                                                                                                                                     \\
$n^{b}_{k,i,j}$  & Number of passengers who can board service $k$ at station $i$ with destination $j$                                                                                                                                  \\
$n^{a}_{k,i}$   & Number of passenger alighting from service $k$ at station $i$                                                                                                                                                     \\
$n_{k,i}$     & Number of in-vehicle passengers after service $k$ serves station $i$                                                                                                                                              \\
$v_{k,i,j}$    & \begin{tabular}[c]{@{}l@{}}Number of remaining passengers who cannot board service $k$ at station $i$ with \\ destination $j$\end{tabular}                                                 \\ \hline
\end{tabular}
}
\end{table}
\\
\textbf{Modelling of passenger demand:}\\
Let $p_{i,j}$ denote the per-second arrival rate of passenger demand during a specific planning time interval. We assume that each demand segment accumulates a two-minute per-second demand before the arrival of the first train at the station. The demand for the second and subsequent trains depends on the headway between trains, as well as the remaining passengers after the previous train operations. Therefore, $w_{k,i,j}$ should be calculated based on the following constraints:
\begin{align}
 &w_{k,i,j} = p_{i,j}.120 & \forall k = 1,i,j \in S^{(up)},i<j\label{48}\\
 &w_{k,i,j} = v_{k-1,i,j} + p_{i,j}.(d_{k,i} - d_{k-1,i}) & \forall k-1, k \in K^{(up)}/\{1\},i,j \in S^{(up)},i<j\label{49}
    \end{align}
Based on Assumption 4, passengers with an O-D pair $i$-$j$ select service $k$ if it serves both the origin station $i$ and the destination station $j$. Variables $w^{b}_{k,i,j}$ and $w_{k,i,j}$ must satisfy:
     \begin{align}
     &w^{b}_{k,i,j} = w_{k,i,j}. x_{k,i}. x_{k,j} & \forall k \in K^{(up)},i,j \in S^{(up)},i<j\label{50}
     \end{align}   
Additionally, following Assumption 4, passenger transfers between services are not allowed, ensuring the following constraints always hold for $w^{b}_{k,i}$ and $w^{b}_{k,i,j}$.
      \begin{align}
     &w^{b}_{k,i} = \sum^{N}_{j=i+1}w^{b}_{k,i,j} & \forall k \in K^{(up)},i\in S^{(up)}\label{51}
     \end{align}
Due to capacity restrictions on trains, some passengers may be unable to board in overcrowded scenarios. The maximum number of passengers permitted on train service $k$ at station $i$, denoted as $n^{b}_{k,i}$, is determined by the following constraint:
      \begin{align}
      &n^{b}_{k,i} = min[C-n_{k,i-1} + n^{a}_{k,i} , w^{b}_{k,i}]  & \forall k \in K^{(up)},i\in S^{(up)}\label{52}
      \end{align}
The number of passengers alighting from train service $k$ at station $i$ is determined by the following constraint:
       \begin{align}
      &n^{a}_{k,i} = \sum^{i-1}_{j=1} n^{b}_{k,j,i}  & \forall k \in K^{(up)},i\in S^{(up)}\label{53}
      \end{align}
Next, we determine $n^{b}_{k,i,j}$, which denotes the count of passengers eligible to board service $k$ at station $i$ for destination $j$. If the number of passengers who want to board service $k$ at station $i$ is equal to the number of passengers who can board service $k$ at station $i$, i.e., $w^{b}_{k,i}$ = $n^{b}_{k,i}$, hence all commuters waiting for service $k$ will be able to board. If $w^{b}_{k,i}$ $>$ $n^{b}_{k,i}$, an appropriate value for $n^{b}_{k,i,j}$ is needed. Passengers board the first available train serving their destination stations if there is space. The constraints listed below establish the relationship between the variables $w^{b}_{k,i}$, $n^{b}_{k,i}$, $n^{b}_{k,i,j}$, and $w^{b}_{k,i,j}$.
  \begin{align}
           &n^{b}_{k,i} - n^{b}_{k,i,j} = w^{b}_{k,i} -w^{b}_{k,i,j}  & \forall k \in K^{(up)},i,j\in S^{(up)},i<j\label{54}
       \end{align}       
After departing from station $i$, the number of passengers on board service $k$, represented as $n_{k,i}$, is calculated as follows:
\begin{align}
&n_{k,i} = n_{k,i-1} - n^{a}_{k,i} + n^{b}_{k,i} & \forall k \in K^{(up)},i\in S^{(up)}\label{55}
  \end{align}
The number of remaining passengers, denoted as $v_{k,i,j}$, who cannot board service $k$ due to capacity restrictions at station $i$ with destination $j$, can be determined as follows:
\begin{align}
&v_{k,i,j} = w_{k,i,j}-n^{b}_{k.i.j} & \forall k \in K^{(up)},i,j\in S^{(up)},i<j\label{eq:56}
\end{align}
As the double-track configuration distinguishes between upstream and downstream services, the demand constraints for services in the downstream direction have the same interpretation as the upstream demand constraints (48-56):
\begin{align}
      &w_{k,i,j} = p_{i,j}.120  &\forall k = 1,i,j \in S^{(dn)},i<j\label{57}\\
      &w_{k,i,j} = v_{k-1,i,j} + p_{i,j}.(d_{k,i} - d_{k-1,i})  &\forall k \in K^{(dn)}/\{1\},i,j \in S^{(dn)},i<j\label{58}\\
      &w^{b}_{k,i,j} = w_{k,i,j}. x_{k,i}. x_{k,j} & \forall k \in K^{(dn)},i,j \in S^{(dn)},i<j\label{59}\\
      &w^{b}_{k,i} = \sum^{N}_{j=i+1}w^{b}_{k,i,j}  &\forall k \in K^{(dn)},i\in S^{(dn)}\label{60}\\
      &n^{b}_{k,i} = min[C-n_{k,i-1} + n^{a}_{k,i} , w^{b}_{k,i}]   &\forall k \in K^{(dn)},i\in S^{(dn)}\label{61}\\
      &c^{a}_{k,i} = \sum^{i-1}_{j=1} c^{b}_{k,j,i}   &\forall k \in K^{(dn)},i\in S^{(dn)}\label{62}\\
      &n^{b}_{k,i} - n^{b}_{k,i,j} = w^{b}_{k,i} -w^{b}_{k,i,j} &\forall k \in K^{(dn)},i,j\in S^{(dn)},i<j\label{63}\\
      &n_{k,i} = n_{k,i-1} - n^{a}_{k,i} + n^{b}_{k,i}  &\forall k \in K^{(dn)},i\in S^{(dn)}\label{64}\\  
      &v_{k,i,j} = w_{k,i,j}-c^{b}_{k,i,j}  &\forall k \in K^{(dn)},i,j \in S^{(dn)},i<j\label{65}
\end{align}
\subsubsection{Peak hours}
Stations often experience increased passenger volumes during specific hours of the day, which leads to overcrowding. 
This congestion requires additional operational strategies beyond the standard planning for non-peak hours. In our initial consideration for regular off-peak hours, we focused on integrating a short-turning strategy along with rolling stock assignment and timetable optimization to enhance operational efficiency and passenger satisfaction.

However, in peak hours, we need additional strategies to reduce overcrowding and minimize passenger waiting times. 
This paper introduces the skip-stop strategy specifically designed for peak hours. 
This strategy involves certain train services skipping specific stations to reduce extra acceleration, deceleration, and dwelling time. It aims to maintain an average maximum speed, improving resource utilization, minimizing congestion, and optimizing service efficiency. Information regarding which train services skip certain stations is communicated well in advance of the peak hour period, and specific trains already in operation adopt this strategy after reaching any turnaround station.

To incorporate skip-stopping, we need to remove some constraints and introduce new constraints. Changes are applied to arrival, departure, and dwell time constraints to accommodate the skip-stopping strategy. Furthermore, some additional constraints are introduced to regulate the skipping behavior.\\
\textbf{New arrival and departure constraints:}\\
The time takes by a train to travel between two consecutive stations, \( i-1 \) and \( i \), can be broken down into three components: the base running time \( Rt_i \), additional time for acceleration \( \tau_a \), and extra time for deceleration \( \tau_d \). When a train skips a station, it can save the additional time needed for acceleration and deceleration ensured by following constraints. This is another reason why trains tend to have a higher average speed when using a skip-stop pattern compared to a standard stop pattern, in addition to the saved dwell time.
\begin{align}
&a_{k,i} - d_{k,i-1} = Rt_{i-1,i} + x_{k,i-1}.R_{a} + x_{k,i}.R_{d}  & \forall k \in K^{(up)}, i-1,i \in S^{(up)} \label{66}\\
&a_{k,i} - d_{k,i-1} = Rt_{i-1,i} + x_{k,i-1}.R_{a} + x_{k,i}.R_{d} & \forall k \in K^{(dn)}, i-1,i \in S^{(dn)} \label{67}
\end{align}
\textbf{New dwell constraints:}\\
In the skip-stop pattern, when service $k$ stops at station $i$, the time spent at the station (dwell time) must be greater than a minimum value to allow for boarding and alighting. If the train skips the station, the dwell time is considered to be 0. Therefore, the dwell time of service $k$  at station $i$, represented by $(d_{k,i} - a_{k,i})$, needs to satisfy following constraints.
\begin{align}
&d_{k,i} - a_{k,i} \geq e_{i}.x_{k,i} & \forall k \in K^{(up)}, i \in S^{(up)} \label{68}\\
&d_{k,i} - a_{k,i} \geq e_{i}.x_{k,i} & \forall k \in K^{(dn)}, i \in S^{(dn)} \label{69}
\end{align}
\textbf{Skip-stopping constraints :}\\
To prevent passenger complaints about excessive station skipping, it is essential to impose limitations on the number of stations skipped in each direction. The maximum allowed number of skipped stations, denoted as $s_{k}$, applies to both upstream and downstream directions for service $k$.
\begin{align}
&\sum_{i\in S^{(up)}} x_{k,i} \geq |N_{k}| - s_{k} & \forall k \in K^{(up)} \label{70}\\
&\sum_{i\in S^{(dn)}} x_{k,i} \geq |N_{k}| - s_{k} & \forall k \in K^{(dn)} \label{71}
\end{align}
Apart from that, we need to remove constraint (\ref{6}) and (\ref{10}), which allow serving all the stations in designated operation zones. All other constraints remain the same to integrate skip-stopping for peak hours. With the help of skip-stopping, we can save time on acceleration, deceleration, and dwelling.

These constraints incorporate the integration of the train timetable, rolling stock assignment, and short-turning strategy for off-peak hours, along with the skip-stop strategy for peak hours. In the following section, we introduce the objective function aimed at minimizing operational costs and passenger waiting time.

\subsection{Objective functions}
\subsubsection{Minimizing operational costs}
The operational costs of a metro line mainly depends on the number of trains required to meet the demand between O-D pairs. To minimize the operational cost of a metro line, it is essential to use the minimum number of trains in operations. This is achieved through turnaround existing trains instead of introducing new once. 
The first objective focuses on maximizing the turnaround operations of trains to minimize operational costs, without considering passenger waiting time.

\begin{align*}
F_{obj,1} = \max~&\sum_{k \in K^{(up)}} \sum_{l\in K^{(dn)}} \sum_{m\in S^{(up)}}y_{k,l,m}+\sum_{k \in K^{(dn)}} \sum_{l\in K^{(up)}} \sum_{m\in S^{(dn)}}y_{k,l,m}&
\end{align*}
\subsubsection{Maximizing quality of service}
The service quality of the metro line is determined by passenger satisfaction, which is directly related to both waiting time on the platform and the travel time to reach the destination. Minimizing passenger waiting time on the platform involves reducing the headway between consecutive trains or increasing the frequency of trains at each station. To achieve this, the objective function considers two components. The first component represents the total travel time of passengers between the first and last stations of the operation zone for the selected potential service. The second component focuses on minimizing the headway between two consecutive trains at each station to reduce passenger waiting time. Both components aim to minimize the travel and waiting times of passengers, ultimately maximizing the quality of service.
\begin{multline}
F_{obj,2} = \min~\sum_{k \in K^{(up)}}\sum_{m \in S^{(up)}/\{N-2,N\}} \sum_{n \in S^{(up)}/\{1,3\}} (d_{k,n}-d_{k-1,m}).(z_{k,m,n}) \\
+ \sum_{k \in K^{(dn)}} \sum_{m \in S^{(dn)}/\{2N-2,2N\}} \sum_{n \in S^{(dn)}/ \{N+1,N+3\}}  (d_{k,n}-d_{k-1,m}).(z_{k,m,n})\\
+ \sum_{k \in K^{(up)}}\sum_{i \in S^{(up)}} (d_{k,i} - d_{k-1,i})+ \sum_{k \in K^{(dn)}}\sum_{i \in S^{(dn)}} (d_{k,i} - d_{k-1,i}) \label{72}
\end{multline}
The linearized version of above objective function and constraints (\ref{50}), (\ref{52}), (\ref{59}), (\ref{61}) are presented in Appendix A.

\section{Mathematical Models}\label{sec6}
In this paper, we present three mathematical models (each for peak and off-peak hours). The first model aims to minimize operational costs without considering passenger waiting time. The second model focuses on minimizing both passenger waiting and travel time. The third model is a multi-objective optimization model that considers both of the aforementioned objectives. For simplicity, we use the following notation: $\tau$ for $(\tau_{k})$, $x$ for $(x_{k,i})$, $y$ for $(y_{k,l,m})$, $z$ for $(z_{k,m,n})$, $h$ for $(h_{k-1,k})$, $a$ for $(a_{k,i})$, and $d$ for $(d_{k,i})$.
\subsection{Model 1(Minimize Operation Costs): }
\noindent\textbf{(i)} For off-peak hours:
\begin{align*}
[\text{$Model-1a$}]:
\max_{\tau,z,x,h,d,a,y}~&F_{obj,1}&
\end{align*}
\hspace{4cm}\text{s.t.}
\hspace{0.5cm}\text{Constraints} (1)-(49), 
 (51), (53)-(58), (60), (62)-(65) and (75)-(85)\\
\textbf{(ii)} For peak hours:
\begin{align*}
[\text{$Model-1b$}]:
\max_{\tau,z,x,h,d,a,y}~&F_{obj,1}&
\end{align*}
\hspace{3cm}\text{s.t.}\hspace{1cm}
\text{Constraints} (1)-(5), (7)-(9), (11)-(12), (17)-(49), (51), (53)-(58), (60), (62)-(71) and (75)-(85)
\subsection{Model 2(Maximize Quality of Service): }
\noindent\textbf{(i)} For off-peak hours:
\begin{align*}
[\text{$Model-2a$}]:
\min_{\tau,z,x,h,d,a,y}~&F^{*}_{obj,2}& 
\end{align*}
\hspace{3cm}\text{s.t.}\hspace{1cm}
\text{Constraints} (1)-(49), (51), (53)-(58), (60), (62)-(65) and (75)-(89)\\
Where $F^{*}_{obj,2}$ is the linearized version of $F_{obj,2}$ (refer Appendix A).\\
\textbf{(ii)} For peak hours:
\begin{align*}
[\text{$Model-2b$}]:
\min_{\tau,z,x,h,d,a,y}~&F^{*}_{obj,2}& 
\end{align*}
\hspace{3cm}\text{s.t.}\hspace{1cm}
\text{Constraints} (1)-(5), (7)-(9), (11)-(12), (17)-(49), (51), (53)-(58), (60), (62)-(71) and (75)-(89)
\subsection{Model 3(Minimize Operation Costs and Maximize Quality of Service):}
\noindent\textbf{(i)} For off peak hours:
\begin{align*}
[\text{$Model-3a$}]:
\max_{\tau,z,x,h,d,a,y}~&F_{obj,1} \& \min_{\tau,z,x,h,d,a,y}F^{*}_{obj,2}& 
\end{align*}
\hspace{3cm}\text{s.t.}\hspace{1cm}
\text{Constraints} (1)-(49), (51), (53)-(58), (60), (62)-(65) and (75)-(89)\\
\textbf{(ii)} For peak hours:
\begin{align*}
[\text{$Model-3b$}]:
\max_{\tau,z,x,h,d,a,y}~&F_{obj,1} \& \min_{\tau,z,x,h,d,a,y}F^{*}_{obj,2}& 
\end{align*}
\hspace{3cm}\text{s.t.}\hspace{1cm}
\text{Constraints} (1)-(5), (7)-(9), (11)-(12), (17)-(49), (51), (53)-(58), (60), (62)-(71) and (75)-(89)
\section{Implementation}\label{sec7}
In this section, we present the implementation results of a small-sized instance to demonstrate the effectiveness and efficiency of the proposed models on a simplified Santiago metro line 1, which consists of sixteen stations (Fig. 8).  All experiments are performed on a workstation with a 3.60 GHz Intel Xeon processor, 64 GB RAM, and six cores. The models are implemented in Python and solved using OS. GUROBI 11.0.0. The computation is configured to terminate if it exceeds four hours (14400 seconds).
\subsection{Numerical experiments on simplified Santiago metro line 1} 
\subsubsection{Set-up}
We have taken data provided in 2023 INFORMS RAS problem solving competition (\cite{RASCompetition2023}) for our computational experiments.
Line 1 of the Santiago metro, also known as the red line, spans the main corridor of the city, covering a total distance of 20 kilometers and serving 27 stations. The layout of Santiago's metro line 1 is depicted in Fig. 7, where turnaround stations are denoted by solid circles at positions (1, 3, 25, 27) for the upstream direction and (28, 30, 52, 54) for the downstream direction. The direction from station 1 to station 27 is termed the upstream direction, while the direction from station 28 to 54 is designated as the downstream direction.

\begin{figure}[h!]
\begin{center}
\includegraphics[width=12cm]{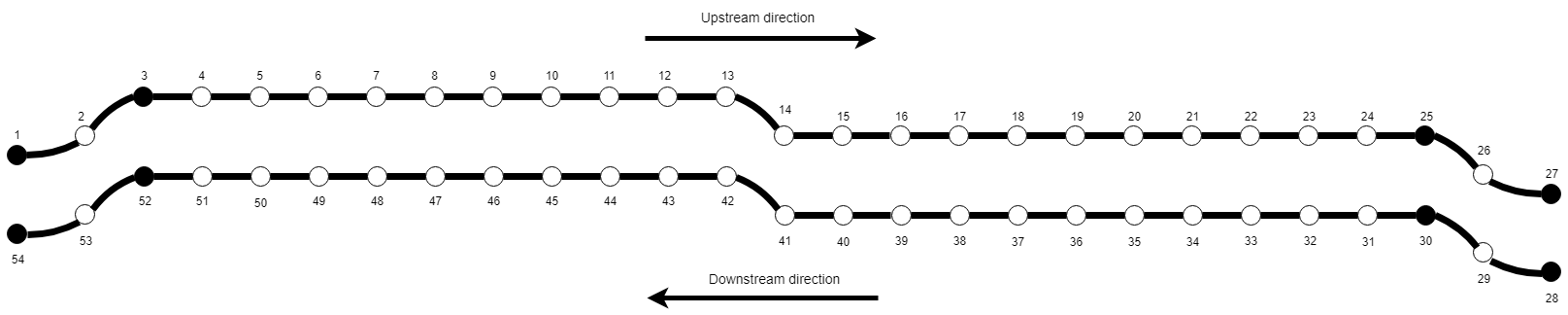}
\end{center}
\caption{Layout of the Santiago metro line 1}
\end{figure}

To assess the proposed models, we conducted small-scale instances on a simplified Santiago metro Line 1, as illustrated in Fig. 8. The bidirectional metro line features sixteen stations, with depots located at both intermediate and terminal stations. Detailed information on dwelling time at stations and running time between adjacent stations is available in Appendix B (refer to Tables 12 and 13). The calculation of the number of potential services involves demand and train capacity calculations, as detailed in subsection 3.2, for both upstream and downstream directions across various time horizons in peak and off-peak hours. The demand matrix for each time horizon is derived from the data provided in the 2023 RAS Problem Solving Competition for forenoon, afternoon, and evening hours. The following Table 5 presents the parameter values used for basic experiments.
\begin{table}[H]
\caption{Parameter values for basic experiments}
\resizebox{\columnwidth}{!}{%
\begin{tabular}{llllll}
\hline
\multicolumn{2}{l}{Headway (seconds)} & \multirow{2}{*}{\begin{tabular}[c]{@{}l@{}}Minimum Turnaround \\ Time (seconds)\end{tabular}} & \multirow{2}{*}{\begin{tabular}[c]{@{}l@{}}Train Capacity\\ (seated and standing)\end{tabular}} & \multicolumn{2}{l}{Load Factor} \\ \cline{1-2} \cline{5-6} 
Maximum & Minimum &  &  & Off-peak hours & Peak hours \\ \hline
360 & 90 & 135 & 250 & 0.8 & 1 \\ \hline
\end{tabular}
}
\end{table}
\begin{figure}[h!]
\begin{center}
\includegraphics[width=12cm]{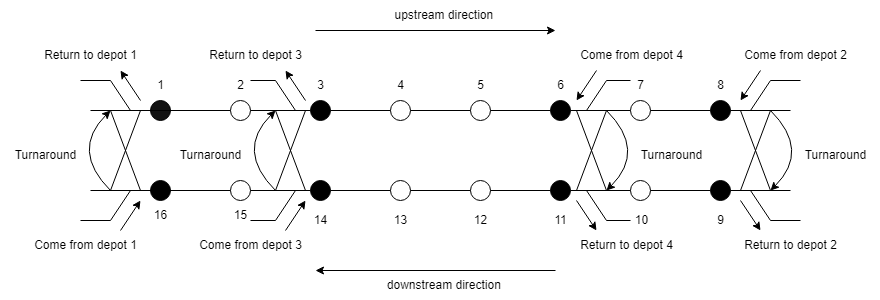}
\end{center}
\caption{A simplified Santiago bidirectional metro line}
\end{figure}
\subsubsection{Basic experiments}
\noindent\textbf{Model 1 (Minimize Operation Costs - for off-peak and peak hours):} In this section, we present the results of Model 1 for both off-peak and peak hours. For off-peak hours, we have taken the demand matrix from the data provided in the 2023 RAS Problem Solving Competition and assumed a 75\% increase in demand for peak hours. The planning time horizon is divided into 15, 30, 45, and 60-minute intervals, spanning from 7:30 am to 8:30 am, 1.00 pm to 2.00 pm and 6.00 pm to 7.00 pm for forenoon, midday and evening, respectively. The passenger arrival rates at stations, represented by $P_{i,j}[t]$, are determined from the origin-destination (O-D) demand matrix. An instance index, denoted as $R-S-T$, outlines the parameter setup where, $R$, $S$, and $T$ signify the number of trains, stations, and time horizon, respectively. For instance, 5-16-15 indicates 5 trains, 16 stations, and 15-minute time frame. Morning, Midday and Evening hours are denoted by M, MD and E, respectively. For peak hour operational planning, a maximum of four stations skipped by a train is considered (i.e. $s_k$ = 4). Table 6 showcases the outcomes of off-peak hour instances with varying train numbers and planning time spans. It shows the objective function value, train allocation across different depots, the count of services chosen in the upstream and downstream directions, along with the CPU time needed to solve each instance.

According to Table 6, to minimize operational costs during morning, midday, and evening hours, an average of 5 trains are needed for both 30 and 60-minute time horizons. Specifically, for the 30-minute time horizon, 7, 4, and 9 turnarounds are optimal during morning, midday, and evening hours, respectively. Likewise, for the 60-minute time horizon, 16, 9, and 21 turnarounds are optimal during morning, midday, and evening hours. 

In Table 7, we present results for peak hours demand where the trains are permitted to implement the skip-stop strategy. We integrate train timetables, short-turning, rolling stock assignment, and the skip-stop strategy into our analysis. During peak hours, the high demand leads to an increase in the number of selected train services in both the up and downstream directions, resulting in an increase in the number of turnarounds. On average, 6 trains and 11 turnarounds are optimal for the morning peak hours, and 12 turnarounds for the evening peak hours with a 30-minute time horizon. Similarly, with a 60-minute time horizon, 22 turnarounds are optimal for the morning peak hours, and 26 for the evening peak hours.

\begin{table}[H]
\caption{Solutions from Model 1}
\resizebox{\columnwidth}{!}{%
\begin{tabular}{llllllllll}
\hline
Instance & \multirow{2}{*}{\begin{tabular}[c]{@{}l@{}}Off-peak hours\\ (M/MD/E)\end{tabular}} & \multicolumn{2}{l}{Number of selected services} & \multicolumn{4}{l}{Number of trains} & \multirow{2}{*}{Objective value} & \multirow{2}{*}{CPU time (sec)} \\ \cline{1-1} \cline{3-8}
R-S-T &  & Upstream & Downstream & Depot 1 & Depot 2 & Depot 3 & Depot 4 &  &  \\ \hline
5-16-30 & M & 6 & 6 & 1 & 3 & 1 & 0 & 7 & 1.6980 \\
5-16-60 & M & 11 & 10 & 2 & 3 & 0 & 0 & 16 & 15.7680 \\
6-16-30 & M & 6 & 6 & 1 & 3 & 1 & 0 & 7 & 1.2850 \\
6-16-60 & M & 11 & 10 & 3 & 1 & 0 & 1 & 16 & 12.1990 \\
7-16-30 & M & 6 & 6 & 2 & 3 & 0 & 0 & 7 & 1.0810 \\
7-16-60 & M & 11 & 10 & 2 & 2 & 1 & 0 & 16 & 16.1580 \\
8-16-30 & M & 6 & 6 & 1 & 3 & 1 & 0 & 7 & 1.4300 \\
8-16-60 & M & 11 & 10 & 1 & 1 & 2 & 1 & 16 & 11.4330 \\
9-16-30 & M & 6 & 6 & 1 & 2 & 1 & 1 & 7 & 1.4770 \\
9-16-60 & M & 11 & 10 & 1 & 1 & 2 & 1 & 16 & 13.4800 \\
10-16-30 & M & 6 & 6 & 1 & 3 & 1 & 0 & 7 & 1.4840 \\
10-16-60 & M & 11 & 10 & 2 & 2 & 0 & 1 & 16 & 5.6020 \\
11-16-30 & M & 6 & 6 & 1 & 3 & 1 & 0 & 7 & 1.5290 \\
11-16-60 & M & 11 & 10 & 2 & 2 & 1 & 0 & 16 & 11.0580 \\
12-16-30 & M & 6 & 6 & 2 & 2 & 1 & 0 & 7 & 1.1460 \\
12-16-60 & M & 11 & 10 & 2 & 3 & 0 & 0 & 16 & 5.9700 \\
13-16-30 & M & 6 & 6 & 2 & 2 & 1 & 0 & 7 & 1.1420 \\
13-16-60 & M & 11 & 10 & 1 & 3 & 1 & 0 & 16 & 5.4120 \\
14-16-30 & M & 6 & 6 & 2 & 2 & 1 & 0 & 7 & 1.1250 \\
14-16-60 & M & 11 & 10 & 1 & 3 & 1 & 0 & 16 & 5.3780 \\
5-16-30 & MD & 5 & 3 & 1 & 1 & 1 & 1 & 4 & 0.1290 \\
5-16-60 & MD & 8 & 6 & 1 & 2 & 2 & 0 & 9 & 1.3700 \\
6-16-30 & MD & 5 & 3 & 2 & 1 & 0 & 1 & 4 & 0.1570 \\
6-16-60 & MD & 8 & 6 & 2 & 2 & 1 & 0 & 9 & 1.4720 \\
7-16-30 & MD & 5 & 3 & 2 & 1 & 0 & 1 & 4 & 0.1420 \\
7-16-60 & MD & 8 & 6 & 2 & 2 & 1 & 0 & 9 & 1.3530 \\
8-16-30 & MD & 5 & 3 & 2 & 1 & 0 & 1 & 4 & 0.1480 \\
8-16-60 & MD & 8 & 6 & 2 & 2 & 1 & 0 & 9 & 1.3720 \\
9-16-30 & MD & 5 & 3 & 2 & 2 & 0 & 0 & 4 & 0.1810 \\
9-16-60 & MD & 8 & 6 & 2 & 2 & 1 & 0 & 9 & 1.3370 \\
10-16-30 & MD & 5 & 3 & 2 & 2 & 0 & 0 & 4 & 0.1440 \\
10-16-60 & MD & 8 & 6 & 2 & 3 & 0 & 0 & 9 & 1.4490 \\
11-16-30 & MD & 5 & 3 & 2 & 2 & 0 & 0 & 4 & 0.1650 \\
11-16-60 & MD & 8 & 6 & 2 & 2 & 1 & 0 & 9 & 1.2300 \\
12-16-30 & MD & 5 & 3 & 2 & 2 & 0 & 0 & 4 & 0.1540 \\
12-16-60 & MD & 8 & 6 & 2 & 3 & 0 & 0 & 9 & 1.3940 \\
13-16-30 & MD & 5 & 3 & 2 & 2 & 0 & 0 & 4 & 0.1510 \\
13-16-60 & MD & 8 & 6 & 2 & 2 & 1 & 0 & 9 & 1.5110 \\
14-16-30 & MD & 5 & 3 & 2 & 2 & 0 & 0 & 4 & 0.1710 \\
14-16-60 & MD & 8 & 6 & 1 & 2 & 1 & 1 & 9 & 0.8490 \\
5-16-30 & E & 6 & 8 & 2 & 2 & 0 & 1 & 9 & 4.7330 \\
5-16-60 & E & 12 & 14 & 2 & 2 & 1 & 0 & 21 & 246.0260 \\
6-16-30 & E & 6 & 8 & 2 & 3 & 0 & 0 & 9 & 4.1230 \\
6-16-60 & E & 12 & 14 & 2 & 3 & 0 & 0 & 21 & 116.2360 \\
7-16-30 & E & 6 & 8 & 1 & 2 & 1 & 1 & 9 & 2.6480 \\
7-16-60 & E & 12 & 14 & 1 & 3 & 1 & 0 & 21 & 78.7230 \\
8-16-30 & E & 6 & 8 & 1 & 2 & 1 & 1 & 9 & 4.5560 \\
8-16-60 & E & 12 & 14 & 3 & 2 & 0 & 0 & 21 & 127.0760 \\
9-16-30 & E & 6 & 8 & 2 & 3 & 0 & 0 & 9 & 1.5440 \\
9-16-60 & E & 12 & 14 & 1 & 3 & 1 & 0 & 21 & 268.6120 \\
10-16-30 & E & 6 & 8 & 2 & 3 & 0 & 0 & 9 & 1.5450 \\
10-16-60 & E & 12 & 14 & 2 & 3 & 0 & 0 & 21 & 86.9540 \\
11-16-30 & E & 6 & 8 & 2 & 2 & 0 & 1 & 9 & 3.5000 \\
11-16-60 & E & 12 & 14 & 3 & 2 & 0 & 0 & 21 & 71.6660 \\
12-16-30 & E & 6 & 8 & 2 & 2 & 0 & 1 & 9 & 4.0130 \\
12-16-60 & E & 12 & 14 & 2 & 3 & 0 & 0 & 21 & 1648.9790 \\
13-16-30 & E & 6 & 8 & 1 & 2 & 1 & 1 & 9 & 3.0730 \\
13-16-60 & E & 12 & 14 & 2 & 3 & 0 & 0 & 21 & 99.3830 \\
14-16-30 & E & 6 & 8 & 2 & 3 & 0 & 0 & 9 & 3.3240 \\
14-16-60 & E & 12 & 14 & 1 & 3 & 1 & 0 & 21 & 314.8910 \\ \hline
\end{tabular}
}
\end{table}
\begin{table}[H]
\caption{Results from Model 1 }
\resizebox{\columnwidth}{!}{%
\begin{tabular}{llllllllll}
\hline
\multicolumn{1}{c}{Instance} & \multirow{2}{*}{\begin{tabular}[c]{@{}l@{}}Peak hours\\ (M/E)\end{tabular}} & \multicolumn{2}{c}{Number of selected services} & \multicolumn{4}{c}{Number of trains} & \multicolumn{1}{c}{\multirow{2}{*}{Objective value}} & \multicolumn{1}{c}{\multirow{2}{*}{CPU time (sec)}} \\ \cline{1-1} \cline{3-8}
\multicolumn{1}{c}{R-S-T} &  & \multicolumn{1}{c}{Upstream} & \multicolumn{1}{c}{Downstream} & \multicolumn{1}{c}{Depot 1} & \multicolumn{1}{c}{Depot 2} & \multicolumn{1}{c}{Depot 3} & \multicolumn{1}{c}{Depot 4} & \multicolumn{1}{c}{} & \multicolumn{1}{c}{} \\ \hline
5-16-30 & M & - & - & - & - & - & - & - & 14400.0720 \\
5-16-60 & M & - & - & - & - & - & - & - & 14400.0720 \\
6-16-30 & M & 9 & 8 & 2 & 1 & 2 & 1 & 11 & 6582.0290 \\
6-16-60 & M & - & - & - & - & - & - & - & 14400.1630 \\
7-16-30 & M & 9 & 8 & 3 & 3 & 0 & 0 & 11 & 12796.0450 \\
7-16-60 & M & 15 & 14 & 5 & 2 & 0 & 0 & 22 &  14400.1390\\
8-16-30 & M & 9 & 8 & 3 & 1 & 1 & 1 & 11 & 6287.9020 \\
8-16-60 & M & 15 & 14 & 4 & 3 & 0 &  0 & 22 &  14400.2610\\
9-16-30 & M & 9 & 8 & 2 & 2 & 1 & 1 & 11 & 11213.3610 \\
9-16-60 & M & 15 & 14 & 2 & 5 & 0 & 0 & 22 &  14400.0990\\
10-16-30 & M & 9 & 8 & 3 & 3 & 0 & 0 & 11 & 11974.2880 \\
10-16-60 & M & 15 & 14 & 5 & 2 & 0 & 0 & 22 &  14400.2960\\
11-16-30 & M & 9 & 8 & 3 & 2 & 1 & 0 & 11 & 4004.0820 \\
11-16-60 & M & 15 & 14 & 3 & 2 & 1 & 0 & 23 &  14400.1700\\
12-16-30 & M & 9 & 8 & 3 & 3 & 0 & 0 & 11 & 9620.0600 \\
12-16-60 & M & 15 & 14 & 4 & 2 & 0 & 1 & 22 &  14400.1850\\
13-16-30 & M & 9 & 8 & 4 & 2 & 0 & 0 & 11 & 7498.0510 \\
13-16-60 & M & 15 & 14 & 2 & 7 & 0 &  0 & 22 &  14400.0810\\
14-16-30 & M & 9 & 8 & 3 & 3 & 0 & 0 & 11 & 8624.3040 \\
14-16-60 & M & 15 & 14 & 3 & 4 & 0 & 0 & 22 & 14400.1160 \\
5-16-30 & E & - & - & - & - & - & - & Infeasible & 1056.7170 \\
5-16-60 & E & - & - & - & - & - & - & - &  14400.1160\\
6-16-30 & E & 8 & 10 & 2 & 4 & 0 & 0 & 12 &10078.0880  \\
6-16-60 & E & - & - & - & - & - & - & - &  14400.1390\\
7-16-30 & E & 8 & 10 & 2 & 4 & 0 & 0 & 12 &  4858.5530\\
7-16-60 & E & - & - & - & - & - & - & - &  14400.1130\\
8-16-30 & E & 8 & 10 & 2 & 4 & 0 & 0 & 12 & 2369.2910\\
8-16-60 & E & - & - & - & - & - & - & - &  14400.1720\\
9-16-30 & E & 8 & 10 & 2 & 4 & 0 & 0 & 12 &  2448.5560\\
9-16-60 & E & - & - & - & - & - & - & - &   14400.1130\\
10-16-30 & E & 8 & 10 & 2 & 3 & 0 & 1 & 12 & 2756.2780 \\
10-16-60 & E & - & - & - & - & - & -& - &  14400.1390\\
11-16-30 & E & 8 & 10 & 2 & 4 & 0 & 0 & 12 & 1313.8390 \\
11-16-60 & E & 16 & 19 & 4 & 4 & 1 & 0 & 26 &  14400.1720\\
12-16-30 & E & 8 & 10 & 2 & 4 & 0 & 0 & 12 & 2347.0890 \\
12-16-60 & E & 16 & 19 & 2 & 7 & 0 & 0 & 26 &  14400.1390\\
13-16-30 & E & 8 & 10 & 2 & 4 & 0 & 0 & 12 & 1213.1960 \\
13-16-60 & E & 16 & 18 & 3 & 6 & 0 & 0 & 22 &14400.0970  \\
14-16-30 & E & 8 & 10 & 2 & 4 & 0 & 0 & 12 & 1234.8770 \\
14-16-60 & E & 16 & 19 & 2 & 4 & 0 & 2 & 22 & 14400.1050 \\ \hline
\multicolumn{10}{l}{Note: -, time limit exceeded but not found sub-optimal solutions; Infeasible, no feasible solution is found for specific instance }
\end{tabular}
}
\end{table}
\noindent\textbf{Model 2 (Maximize Service Quality - for off-peak and peak hours):} Tables 8 and 9 presents the results for off-peak and peak hours, aiming to minimize passengers waiting and travel time. According to Table 8, to minimize passenger waiting and travel time, an average of 6, 3, and 7 trains are required during morning, midday, and evening off-peak hours for 30-minute time horizons, respectively. Similarly, for a 60-minute time horizon, 10, 7, and 13 trains are required during morning, midday, and evening off-peak hours.

During peak hours, an increasing number of services are necessary to minimize passenger waiting time. According to Table 9, all available trains from depots are utilized for operations. Additionally, the problem exhibits significant time complexity in finding solutions, as almost all instances reach the maximum solution time limit of four hours.

\begin{table}[H]
\caption{Solutions from Model 2}
\resizebox{\columnwidth}{!}{%
\begin{tabular}{llllllllll}
\hline
Instance & \multirow{2}{*}{\begin{tabular}[c]{@{}l@{}}Off-peak hours\\ (M/MD/E)\end{tabular}} & \multicolumn{2}{l}{Number of selected services} & \multicolumn{4}{l}{Number of trains} & \multirow{2}{*}{Objective value} & \multirow{2}{*}{CPU time (sec)} \\ \cline{1-1} \cline{3-8}
R-S-T &  & Upstream & Downstream & Depot 1 & Depot 2 & Depot 3 & Depot 4 &  &  \\ \hline
5-16-30 & M & 3 & 3 & 3 & 2 & 0 & 0 & 6559.8222 & 0.3140 \\
5-16-60 & M & 5 & 5 & 2 & 3 & 0 & 0 & 26478.7554 & 82.3280 \\
6-16-30 & M & 3 & 3 & 3 & 3 & 0 & 0 & 6559.8222 & 0.3140 \\
6-16-60 & M & 5 & 5 & 2 & 4 & 0 & 0 & 20985.8962 & 31.9040 \\
7-16-30 & M & 3 & 3 & 3 & 3 & 0 & 0 & 6559.8222 & 0.3090 \\
7-16-60 & M & 5 & 5 & 2 & 5 & 0 & 0 & 16799.4666 & 6.8220 \\
8-16-30 & M & 3 & 3 & 3 & 3 & 0 & 0 & 6559.8222 & 0.3070 \\
8-16-60 & M & 5 & 5 & 3 & 5 & 0 & 0 & 16079.4666 & 6.3820 \\
9-16-30 & M & 3 & 3 & 3 & 3 & 0 & 0 & 6559.8222 & 0.2740 \\
9-16-60 & M & 5 & 5 & 4 & 5 & 0 & 0 & 15359.4666 & 3.7550 \\
10-16-30 & M & 3 & 3 & 3 & 3 & 0 & 0 & 6559.8222 & 0.2980 \\
10-16-60 & M & 5 & 5 & 5 & 5 & 0 & 0 & 12613.0370 & 0.6210 \\
11-16-30 & M & 3 & 3 & 3 & 3 & 0 & 0 & 6559.8222 & 0.3150 \\
11-16-60 & M & 5 & 5 & 5 & 5 & 0 & 0 & 12613.0370 & 0.6540 \\
12-16-30 & M & 3 & 3 & 3 & 3 & 0 & 0 & 6559.8222 & 0.3380 \\
12-16-60 & M & 5 & 5 & 5 & 5 & 0 & 0 & 12613.0370 & 0.5740 \\
13-16-30 & M & 3 & 3 & 3 & 3 & 0 & 0 & 6559.8222 & 0.3450 \\
13-16-60 & M & 5 & 5 & 5 & 5 & 0 & 0 & 12613.0370 & 0.6920 \\
14-16-30 & M & 3 & 3 & 3 & 3 & 0 & 0 & 6559.8222 & 0.2740 \\
14-16-60 & M & 5 & 5 & 5 & 5 & 0 & 0 & 12613.0370 & 0.6600 \\
5-16-30 & MD & 2 & 1 & 2 & 1 & 0 & 0 & 3999.9111 & 0.1350 \\
5-16-60 & MD & 4 & 3 & 2 & 3 & 0 & 0 & 12799.5555 & 1.5000 \\
6-16-30 & MD & 2 & 1 & 2 & 1 & 0 & 0 & 3999.9111 & 0.1610 \\
6-16-60 & MD & 4 & 3 & 3 & 3 & 0 & 0 & 12079.5555 & 1.6060 \\
7-16-30 & MD & 2 & 1 & 2 & 1 & 0 & 0 & 3999.9111 & 0.1550 \\
7-16-60 & MD & 4 & 3 & 4 & 3 & 0 & 0 & 7893.1259 & 0.2960 \\
8-16-30 & MD & 2 & 1 & 2 & 1 & 0 & 0 & 3999.9111 & 0.1690 \\
8-16-60 & MD & 4 & 3 & 4 & 3 & 0 & 0 & 7893.1259 & 0.3100 \\
9-16-30 & MD & 2 & 1 & 2 & 1 & 0 & 0 & 3999.9111 & 0.1260 \\
9-16-60 & MD & 4 & 3 & 4 & 3 & 0 & 0 & 7893.1259 & 0.4000 \\
10-16-30 & MD & 2 & 1 & 2 & 1 & 0 & 0 & 3999.9111 & 0.1720 \\
10-16-60 & MD & 4 & 3 & 4 & 3 & 0 & 0 & 7893.1259 & 0.3730 \\
11-16-30 & MD & 2 & 1 & 2 & 1 & 0 & 0 & 3999.9111 & 0.1650 \\
11-16-60 & MD & 4 & 3 & 4 & 3 & 0 & 0 & 7893.1259 & 0.3580 \\
12-16-30 & MD & 2 & 1 & 2 & 1 & 0 & 0 & 3999.9111 & 0.1500 \\
12-16-60 & MD & 4 & 3 & 4 & 3 & 0 & 0 & 7893.1259 & 0.3460 \\
13-16-30 & MD & 2 & 1 & 2 & 1 & 0 & 0 & 3999.9111 & 0.1350 \\
13-16-60 & MD & 4 & 3 & 4 & 3 & 0 & 0 & 7893.1259 & 0.3500 \\
14-16-30 & MD & 2 & 1 & 2 & 1 & 0 & 0 & 3999.9111 & 0.1400 \\
14-16-60 & MD & 4 & 3 & 4 & 3 & 0 & 0 & 7893.1259 & 0.3240 \\
5-16-30 & E & 3 & 4 & 3 & 2 & 0 & 0 & 12799.5555 & 1.8830 \\
5-16-60 & E & 6 & 7 & 2 & 3 & 0 & 0 & 35385.0962 & 1885.5270 \\
6-16-30 & E & 3 & 4 & 3 & 3 & 0 & 0 & 12079.5555 & 1.3600 \\
6-16-60 & E & 6 & 7 & 3 & 3 & 0 & 0 & 30478.6665 & 2813.7110 \\
7-16-30 & E & 3 & 4 & 3 & 4 & 0 & 0 & 7893.1259 & 0.3000 \\
7-16-60 & E & 6 & 7 & 4 & 3 & 0 & 0 & 24852.2369 & 72.1590 \\
8-16-30 & E & 3 & 4 & 3 & 4 & 0 & 0 & 7893.1259 & 0.3160 \\
8-16-60 & E & 6 & 7 & 6 & 2 & 0 & 0 & 20799.3777 & 19.1770 \\
9-16-30 & E & 3 & 4 & 3 & 4 & 0 & 0 & 7893.1259 & 0.3400 \\
9-16-60 & E & 6 & 7 & 6 & 3 & 0 & 0 & 20079.3777 & 20.1960 \\
10-16-30 & E & 3 & 4 & 3 & 4 & 0 & 0 & 7893.1259 & 0.3390 \\
10-16-60 & E & 6 & 7 & 6 & 4 & 0 & 0 & 19359.3777 & 17.6360 \\
11-16-30 & E & 3 & 4 & 3 & 4 & 0 & 0 & 7893.1259 & 0.3480 \\
11-16-60 & E & 6 & 7 & 6 & 5 & 0 & 0 & 18639.3777 & 9.0080 \\
12-16-30 & E & 3 & 4 & 3 & 4 & 0 & 0 & 7893.1259 & 0.3440 \\
12-16-60 & E & 6 & 7 & 6 & 6 & 0 & 0 & 17919.3777 & 2.4640 \\
13-16-30 & E & 3 & 4 & 3 & 4 & 0 & 0 & 7893.1259 & 0.3020 \\
13-16-60 & E & 6 & 7 & 6 & 7 & 0 & 0 & 15892.9481 & 0.7790 \\
14-16-30 & E & 3 & 4 & 3 & 4 & 0 & 0 & 7893.1259 & 0.3550 \\
14-16-60 & E & 6 & 7 & 6 & 7 & 0 & 0 & 15892.9481 & 0.7220 \\ \hline
\end{tabular}
}
\end{table}

\begin{table}[H]
\caption{Results from Model 2}
\resizebox{\columnwidth}{!}{%
\begin{tabular}{llllllllll}
\hline
\multicolumn{1}{c}{Instance} & \multirow{2}{*}{\begin{tabular}[c]{@{}l@{}}Peak hours\\ (M/E)\end{tabular}} & \multicolumn{2}{c}{Number of selected services} & \multicolumn{4}{c}{Number of trains} & \multicolumn{1}{c}{\multirow{2}{*}{Objective value}} & \multicolumn{1}{c}{\multirow{2}{*}{CPU time (sec)}} \\ \cline{1-1} \cline{3-8}
\multicolumn{1}{c}{R-S-T} &  & \multicolumn{1}{c}{Upstream} & \multicolumn{1}{c}{Downstream} & \multicolumn{1}{c}{Depot 1} & \multicolumn{1}{c}{Depot 2} & \multicolumn{1}{c}{Depot 3} & \multicolumn{1}{c}{Depot 4} & \multicolumn{1}{c}{} & \multicolumn{1}{c}{} \\ \hline
5-16-30 & M & - & - & - & - & - & - & Infeasible & 6995.9270 \\
5-16-60 & M & - & - & - & - & - & - & - &  14400.2360\\
6-16-30 & M & 9 & 8 & 2 & 2 & 1 & 1 & 26442.0137 & 14400.3540 \\
6-16-60 & M & - & - & - & - & - & - & - & 14400.1650 \\
7-16-30 & M & 9 & 8 & 2 & 2 & 1 & 2 & 25076.9235 & 14400.2240 \\
7-16-60 & M & 15 & 13 & 2 & 2 & 1 & 2 &  51011.1967 &  14400.2360\\
8-16-30 & M & 9 & 8 & 1 & 3 & 1 & 3 & 21319.9468 & 14400.4200 \\
8-16-60 & M & 15 & 14 & 2 & 3 & 1 &  2 & 44207.7983 &  14400.1530\\
9-16-30 & M & 9 & 8 & 3 & 3 & 1 & 2 & 20416.2303 & 14400.8710 \\
9-16-60 & M & 15 & 14 & 3 & 3 & 2 & 1 & 38315.4808 &  14400.2530\\
10-16-30 & M & 9 & 8 & 3 & 4 & 1 & 2 & 19242.9845 & 14400.1350 \\
10-16-60 & M & 15 & 14 & 3 & 3 & 1 & 3 & 38056.2216 & 14400.2800 \\
11-16-30 & M & 9 & 8 & 4 & 4 & 1 & 2 & 17884.0451 & 2049.1560 \\
11-16-60 & M & 15 & 14 & 3 & 4 & 1 & 3 & 35346.7594 &  14400.3230\\
12-16-30 & M & 9 & 8 & 5 & 4 & 1 & 2 & 17506.3958 & 680.3000 \\
12-16-60 & M & 15 & 13 & 3 & 5 & 2 & 2 & 33325.8812 &  14400.3400\\
13-16-30 & M & 9 & 8 & 5 & 4 & 2 & 2 & 17154.9149 & 427.8890 \\
13-16-60 & M & 15 & 14 & 4 & 5 & 2 & 2 & 31580.0201 &  14400.2450\\
14-16-30 & M & 9 & 8 & 4 & 4 & 3 & 3 & 17143.8078 & 296.6710 \\
14-16-60 & M & 15 & 14 & 4 & 4 & 3 & 3 & 30264.8275 & 14400.2780 \\
5-16-30 & E & - & - & - & - & - & - & Infeasible & 8350.4120 \\
5-16-60 & E & - & - & - & - & - & - & - &  14400.3400\\
6-16-30 & E & 8 & 10 & 2 & 2 & 1 & 1 & 29234.9672 & 14400.1450 \\
6-16-60 & E & - & - & - & - & - & - & - &  14400.2780 \\
7-16-30 & E & 8 & 10 & 2 & 2 & 1 & 2 & 26850.7017 & 14400.2970 \\
7-16-60 & E & - & - & - & - & - & - & - &  14400.1450\\
8-16-30 & E & 8 & 10 & 3 & 2 & 2 & 1 & 22810.6491 & 14400.5270 \\
8-16-60 & E & - & - & - & - & - & - & - &  14400.1370\\
9-16-30 & E & 8 & 10 & 4 & 1 & 3 & 1 & 21767.0992 & 14400.4940 \\
9-16-60 & E & - & - & - & - & - & - & - & 14400.4770 \\
10-16-30 & E & 8 & 10 & 4 & 3 & 2 & 1 & 20161.1598 & 14400.2470 \\
10-16-60 & E & 16 & 19 & 4 & 3 & 1 & 2 & 45220.6531 & 14400.4860 \\
11-16-30 & E & 8 & 10 & 4 & 4 & 2 & 1 & 19179.2309 & 7874.9110 \\
11-16-60 & E & 16 & 19 & 3 & 5 & 2 & 1 & 45650.1106 & 14400.3150 \\
12-16-30 & E & 8 & 10 & 4 & 4 & 2 & 2 & 18692.2765 & 2107.4980 \\
12-16-60 & E & 16 & 19 & 4 & 4 & 2 & 2 & 41911.2704 & 14400.3140 \\
13-16-30 & E & 8 & 10 & 4 & 5 & 2 & 2 & 18239.0296 & 298.1870 \\
13-16-60 & E & 16 & 19 & 6 & 3 & 3 & 1 & 38900.9071 &  14400.2440\\
14-16-30 & E & 8 & 10 & 4 & 5 & 2 & 3 & 18239.0296 & 243.6190 \\
14-16-60 & E & 16 & 19 & 5 & 4 & 3 & 2 & 38000.6939 & 14400.3680 \\ \hline
\multicolumn{10}{l}{Note: -, time limit exceeded but not found sub-optimal solutions; Infeasible, no feasible solution is found for specific instance }

\end{tabular}
}
\end{table}
\noindent\textbf{Model 3 (Minimize Operation Costs and Maximize Service Quality - for off-peak and peak hours):} Unlike single-objective optimization models, multi-objective optimization models do not have a single solution that represents the best values for all considered objectives. Specifically, the mathematical formulation of Model 3 involves two objective functions that conflict with each other. Improving one objective function will result in a worsening of the value of the other objective function. The epsilon constraint approach is used to solve multi-objective optimization models. This method was initially introduced by \cite{haimes1971bicriterion}, and later discussed extensively by \cite{chankong2008multiobjective}. In this approach, one objective function is chosen as the primary objective, and remaining objectives are converted into constraints. We use epsilon constraint approach to solve the multi-objective optimization model 3. We first solve for objective 1 while taking all constraints into account, which yields values for the decision variables. These values are subsequently substituted into the expression for objective function 2.
\begin{align}
\max_{\tau,z,x,h,d,a,y}~&F_{obj,1}   &\label{73}
\end{align}
\hspace{3cm}\text{s.t.}\hspace{1cm}
\text{Constraints} (1)-(49), 
 (51), (53)-(58), (60), (62)-(65) and (75)-(85)\\
The resulting cut (74) serves as an epsilon constraint, which is generated and added with the existing constraints in each iteration. This iterative process continues with sequential cut generation until the problem becomes infeasible.
\begin{align}
 &F^{*}_{obj,2} \leq F^{**}_{obj,2} - \epsilon & \label{74}
\end{align}
where $F^{**}_{obj,2}$ represents the value of objective function 2 by substituting the obtained decision variables.
\begin{figure}[h!]
\begin{center}
\includegraphics[width=15cm]{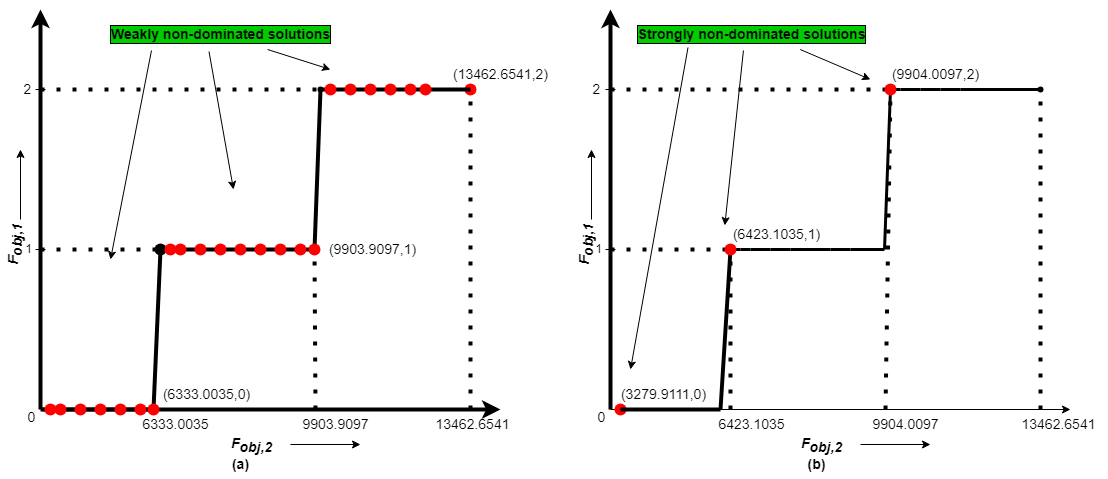}
\end{center}
\caption{(a) Weakly non-dominated solutions and (b) Strongly non-dominated solutions}
\end{figure}
\begin{table}[H]
\caption{Results from Model 3}
\resizebox{\columnwidth}{!}{%
\begin{tabular}{lllllllllll}
\hline
Instance & \multirow{2}{*}{\begin{tabular}[c]{@{}l@{}}Off-peak hours\\ (M/MD/E)\end{tabular}} & \multicolumn{2}{l}{Non-dominated Solutions} & \multicolumn{2}{l}{Number of selected services} & \multicolumn{4}{l}{Number of trains} & \multirow{2}{*}{CPU time (sec)} \\ \cline{1-1} \cline{3-10}
R-S-T &  & Objective 1 & Objective 2 & Upstream & Downstream & Depot 1 & Depot 2 & Depot 3 & Depot 4 &  \\ \hline
5-16-30 & M & 1 & 11466.2518 & 3 & 3 & 3 & 2 & 0 & 0 & 1.1830 \\
5-16-30 & M & 2 & 12448.9567 & 3 & 4 & 3 & 1 & 0 & 1 & 1.6550 \\
5-16-30 & M & 3 & 13436.6616 & 4 & 4 & 1 & 3 & 0 & 1 & 2.1490 \\
5-16-30 & M & 4 & 14602.7195 & 5 & 4 & 1 & 3 & 0 & 1 & 2.7080 \\
5-16-30 & M & 5 & 18611.4171 & 6 & 4 & 1 & 3 & 1 & 0 & 8.7290 \\
5-16-30 & M & 6 & 20250.8379 & 6 & 5 & 1 & 2 & 1 & 1 & 3.8710 \\
5-16-30 & M & 7 & 22202.2586 & 6 & 6 & 1 & 2 & 1 & 1 & 2.8180 \\
10-16-60 & M & 0 & 12613.0370 & 5 & 5 & 5 & 5 & 0 & 0 & 1.5220 \\
10-16-60 & M & 1 & 15359.4666 & 5 & 5 & 4 & 5 & 0 & 0 & 12.3640 \\
10-16-60 & M & 2 & 16079.4666 & 5 & 5 & 3 & 5 & 0 & 0 & 18.9880 \\
10-16-60 & M & 3 & 16799.4666 & 5 & 5 & 2 & 5 & 0 & 0 & 36.7520 \\
10-16-60 & M & 4 & 17787.1715 & 6 & 5 & 1 & 5 & 1 & 0 & 73.7070 \\
10-16-60 & M & 5 & 18953.2294 & 7 & 5 & 1 & 5 & 1 & 0 & 82.4760 \\
10-16-60 & M & 6 & 20095.7826 & 8 & 6 & 1 & 5 & 2 & 0 & 120.1990 \\
10-16-60 & M & 7 & 21007.9912 & 8 & 7 & 1 & 5 & 1 & 1 & 99.0560 \\
10-16-60 & M & 8 & 22174.0491 & 9 & 7 & 1 & 5 & 1 & 1 & 124.2600 \\
10-16-60 & M & 9 & 23494.9553 & 10 & 8 & 2 & 5 & 1 & 1 & 457.4400 \\
10-16-60 & M & 10 & 24731.5095 & 11 & 8 & 2 & 5 & 1 & 1 & 797.2870 \\
10-16-60 & M & 11 & 26680.1913 & 10 & 10 & 3 & 3 & 2 & 1 & 4737.9380 \\
10-16-60 & M & 12 & 28287.9246 & 11 & 10 & 2 & 3 & 2 & 2 & 2052.1700 \\
10-16-60 & M & 13 & 31228.0648 & 11 & 10 & 2 & 3 & 1 & 2 & 12193.4230 \\
10-16-60 & M & 14 & 33922.5764 & 11 & 10 & 2 & 2 & 2 & 1 & 5719.1250 \\
10-16-60 & M & 15 & 38487.6641 & 11 & 10 & 2 & 2 & 1 & 1 & 1388.1300 \\
10-16-60 & M & 16 & 45906.2397 & 11 & 10 & 1 & 2 & 1 & 1 & 52.7470 \\
5-16-30 & MD & 0 & 3999.9111 & 2 & 1 & 2 & 1 & 0 & 0 & 0.1290 \\
5-16-30 & MD & 1 & 6423.1035 & 3 & 2 & 2 & 1 & 0 & 1 & 0.2250 \\
5-16-30 & MD & 2 & 9113.5133 & 4 & 3 & 2 & 1 & 0 & 2 & 0.2120 \\
5-16-30 & MD & 3 & 11773.8663 & 4 & 3 & 1 & 1 & 0 & 2 & 0.2910 \\
5-16-30 & MD & 4 & 14695.8580 & 5 & 3 & 1 & 1 & 1 & 1 & 0.1730 \\
7-16-60 & MD & 0 & 7893.1259 & 4 & 3 & 4 & 3 & 0 & 0 & 0.4330 \\
7-16-60 & MD & 1 & 12079.5555 & 4 & 3 & 3 & 3 & 0 & 0 & 2.5810 \\
7-16-60 & MD & 2 & 12799.5555 & 4 & 3 & 2 & 3 & 0 & 0 & 2.7050 \\
7-16-60 & MD & 3 & 13787.2604 & 5 & 3 & 1 & 3 & 1 & 0 & 4.3320 \\
7-16-60 & MD & 4 & 14769.9653 & 5 & 4 & 1 & 3 & 0 & 1 & 3.9110 \\
7-16-60 & MD & 5 & 16740.3752 & 6 & 5 & 1 & 3 & 0 & 2 & 10.0780 \\
7-16-60 & MD & 6 & 17798.5765 & 7 & 5 & 1 & 3 & 1 & 1 & 8.9010 \\
7-16-60 & MD & 7 & 19048.9863 & 8 & 6 & 1 & 3 & 1 & 2 & 5.4470 \\
7-16-60 & MD & 8 & 23688.1985 & 8 & 6 & 1 & 2 & 1 & 2 & 10.9020 \\
7-16-60 & MD & 9 & 27136.1265 & 8 & 6 & 2 & 1 & 1 & 1 & 3.6860 \\
7-16-30 & E & 0 & 7893.1259 & 3 & 4 & 3 & 4 & 0 & 0 & 0.2820 \\
7-16-30 & E & 1 & 12079.5555 & 3 & 4 & 3 & 3 & 0 & 0 & 2.4210 \\
7-16-30 & E & 2 & 12799.5555 & 3 & 4 & 3 & 2 & 0 & 0 & 3.3600 \\
7-16-30 & E & 3 & 13782.2604 & 3 & 5 & 3 & 1 & 0 & 1 & 2.5160 \\
7-16-30 & E & 4 & 15398.9828 & 4 & 6 & 3 & 2 & 1 & 0 & 4.1930 \\
7-16-30 & E & 5 & 16735.8627 & 5 & 7 & 3 & 2 & 1 & 1 & 9.3410 \\
7-16-30 & E & 6 & 17455.8627 & 5 & 7 & 3 & 1 & 1 & 1 & 5.2780 \\
7-16-30 & E & 7 & 18973.4900 & 6 & 8 & 3 & 1 & 2 & 1 & 5.4280 \\
7-16-30 & E & 8 & 23479.1318 & 6 & 8 & 2 & 1 & 2 & 1 & 9.8980 \\
7-16-30 & E & 9 & 28856.8455 & 6 & 8 & 1 & 2 & 1 & 1 & 3.3790 \\
13-16-60 & E & 0 & 15892.9481 & 6 & 7 & 6 & 7 & 0 & 0 & 0.9700 \\
13-16-60 & E & 1 & 17919.3777 & 6 & 7 & 6 & 6 & 0 & 0 & 5.8940 \\
13-16-60 & E & 2 & 18639.3777 & 6 & 7 & 6 & 5 & 0 & 0 & 31.5450 \\
13-16-60 & E & 3 & 19359.3777 & 6 & 7 & 6 & 4 & 0 & 0 & 45.8780 \\
13-16-60 & E & 4 & 20079.3777 & 6 & 7 & 6 & 3 & 0 & 0 & 66.0280 \\
13-16-60 & E & 5 & 20799.3777 & 6 & 7 & 6 & 2 & 0 & 0 & 49.5970 \\
13-16-60 & E & 6 & 21782.0826 & 6 & 8 & 6 & 1 & 0 & 1 & 97.8460 \\
13-16-60 & E & 7 & 23032.4925 & 7 & 9 & 6 & 2 & 0 & 1 & 195.4250 \\
13-16-60 & E & 8 & 24015.1974 & 7 & 10 & 6 & 1 & 0 & 2 & 173.1930 \\
13-16-60 & E & 9 & 25265.6073 & 8 & 11 & 6 & 2 & 0 & 2 & 852.4580 \\
13-16-60 & E & 10 & 26253.3122 & 9 & 11 & 6 & 2 & 1 & 1 & 1601.8870 \\
13-16-60 & E & 11 & 27398.2630 & 9 & 12 & 6 & 2 & 1 & 1 & 2679.7320 \\
13-16-60 & E & 12 & 28352.8566 & 10 & 13 & 5 & 3 & 1 & 2 & 1945.7890 \\
13-16-60 & E & 13 & 29603.2665 & 11 & 14 & 5 & 3 & 1 & 3 & 4368.2120 \\
13-16-60 & E & 14 & 30661.4677 & 12 & 14 & 4 & 4 & 2 & 2 & 10000.0290 \\
13-16-60 & E & 15 & 32030.9714 & 12 & 14 & 4 & 3 & 3 & 1 & 12263.2410 \\
13-16-60 & E & 16 & 33994.2503 & 12 & 14 & 3 & 3 & 2 & 2 & 14400.3520 \\
13-16-60 & E & 17 & 37282.4604 & 12 & 14 & 3 & 3 & 2 & 1 & 14400.3150 \\
13-16-60 & E & 18 & 40162.4604 & 12 & 14 & 3 & 2 & 2 & 1 & 14400.3890 \\
13-16-60 & E & 19 & 43361.6725 & 12 & 14 & 2 & 2 & 2 & 1 & 14400.4720 \\
13-16-60 & E & 20 & 49255.3196 & 12 & 14 & 2 & 2 & 1 & 1 &14400.5340  \\
13-16-60 & E & 21 & 58811.7654  & 12 & 14 & 1 & 2 & 1 & 1 & 338.5870 \\ \hline
\end{tabular}
}
\end{table}

\begin{table}[H]
\caption{Results from Model 3}
\resizebox{\textwidth}{!}{%
\begin{tabular}{lllllllllll}
\hline
Instance & \multirow{2}{*}{\begin{tabular}[c]{@{}l@{}}Peak hours\\ (M/E)\end{tabular}} & \multicolumn{2}{l}{Non-dominated Solutions} & \multicolumn{2}{l}{Number of selected services} & \multicolumn{4}{l}{Number of trains} & \multirow{2}{*}{CPU time (sec)} \\ \cline{1-1} \cline{3-10}
R-S-T &  & Objective 1 & Objective 2 & Upstream & Downstream & Depot 1 & Depot 2 & Depot 3 & Depot 4 &  \\ \hline
30-16-14 & M & 3 & 17143.8078 & 9 & 8 & 4 & 4 & 3 & 3 & 692.6670 \\
30-16-14 & M & 4 & 17154.9149 & 9 & 8 & 5 & 4 & 2 & 2 & 647.3090 \\
30-16-14 & M & 5 & 17506.3958 & 9 & 8 & 5 & 4 & 1 & 2 & 3385.7520 \\
30-16-14 & M & 6 & 17884.0451 & 9 & 8 & 4 & 4 & 1 & 2 & 12603.4060 \\
30-16-14 & M & 7 & 19296.7256 & 9 & 8 & 2 & 4 & 1 & 3 & 14400.3470 \\
30-16-14 & M & 8 & 20347.3577 & 9 & 8 & 2 & 3 & 2 & 2 & 14400.5940 \\
30-16-14 & M & 9 & 21726.5343 & 9 & 8 & 2 & 3 & 1 & 2 & 14400.3280 \\
30-16-14 & M & 10 &  25076.9235 & 9 & 8 & 2 & 2 & 1 & 2 & 14400.3730 \\
30-16-14 & M & 11 &  26442.0136 & 9 & 8 & 2 & 2 & 1 & 1 & 14400.7220 \\
60-16-14 & M & 15 & 30264.8275 & 15 & 14 & 4 & 4 & 3 & 3 & 14400.2520 \\
60-16-14 & M & 16 & 31026.3820 & 15 & 14 & 3 & 5 & 2 & 4 & 14400.2850 \\
60-16-14 & M & 17 & 31572.0739 & 15 & 14 & 4 & 4 & 2 & 2 & 14400.2010 \\
60-16-14 & M & 18 & 35215.5859 & 15 & 14 & 4 & 3 & 2 & 2 & 14400.2610 \\
60-16-14 & M & 19 & 36785.3736 & 15 & 14 & 4 & 3 & 1 & 2 & 14400.5260 \\
60-16-14 & M & 20 & 43151.7417 & 15 & 14 & 3 & 3 & 1 & 2 & 14400.5060 \\
60-16-14 & M & 21 & 41821.3162 & 15 & 14 & 2 & 3 & 1 & 2 & 14400.4650 \\
60-16-14 & M & - & - & - & - & - & - & - & - & 14400.3500 \\
30-16-14 & E & 4 & 18239.0296 & 8 & 10 & 4 & 5 & 2 & 3 & 358.9150 \\
30-16-14 & E & 5 & 18239.0296 & 8 & 10 & 4 & 5 & 2 & 2 & 321.2740 \\
30-16-14 & E & 6 & 18692.2765 & 8 & 10 & 4 & 4 & 2 & 2 & 2084.2790 \\
30-16-14 & E & 7 & 19179.2309 & 8 & 10 & 4 & 4 & 2 & 1 & 11333.4640 \\
30-16-14 & E & 8 & 20161.1598 & 8 & 10 & 4 & 3 & 2 & 1 & 14400.8190 \\
30-16-14 & E & 9 & 21343.6729 & 8 & 10 & 4 & 2 & 2 & 1 & 14400.3580 \\
30-16-14 & E & 10 & 22415.1689 & 8 & 10 & 3 & 1 & 3 & 1 & 14402.2620 \\
30-16-14 & E & 11 & 26254.4512 & 8 & 10 & 2 & 2 & 1 & 3 & 14401.2460 \\
30-16-14 & E & 12 & 29234.9672 & 8 & 10 & 2 & 2 & 1 & 1 & 14401.0760 \\
60-16-14 & E & 21 & 41198.7420 & 16 & 19 & 3 & 5 & 2 & 4 & 14400.7560 \\
60-16-14 & E & 22 & 38207.0142 & 16 & 19 & 4 & 3 & 3 & 3 & 14400.4500 \\
60-16-14 & E & 23 & 39431.2230 & 16 & 19 & 3 & 3 & 3 & 3 & 14401.4600 \\
60-16-14 & E & 24 & 45510.5082 &  16 & 19  & 4 & 5 & 1 & 1 & 14400.1520 \\
60-16-14 & E & 25 & 47682.1680 & 16 & 19 & 4 & 3 & 2 & 1 & 14400.1710 \\
60-16-14 & E & 26 & 51102.4540 & 16 & 19 & 3 & 2 & 3 & 1 & 14400.3940 \\
60-16-14 & E & 27 & - & - & - & - & - & - & - & 14400.4200 \\
 \hline
 \multicolumn{10}{l}{Note: -, time limit exceeded but not found sub-optimal solutions }

\end{tabular}
}
\end{table}

Implementing the above approach, we identify non-dominating solutions with respect to both  objectives. By optimizing the first objective and treating the second objective as a constraint, we generated several solutions, with a majority being weakly non-dominated solutions, as depicted in Figure 9(a) (Instance of 5 trains, 16 station, 15 minute time horizon for Off-peak hours Model 3). For smaller instances like (5-16-15), the addition of multiple cuts is necessary, making the process time-consuming.

On the other hand optimizing objective 2 using objective 1 as an epsilon constraint requires fewer cuts, significantly reducing the solution time. Additionally, we obtained strongly non-dominated solutions for both maximizing $F_{obj,1}$ and minimizing $F_{obj,2}$ (as shown in Figure 9(b)). Notably, since the expression of objective 1 exclusively involves binary variables, we set the epsilon value to 1. The same experimental setup as model 1 is employed for model 3. The results of model 3 are presented in Table 10 (off-peak hours) and Table 11 (for peak hours), showcasing the non-dominated solutions for each specific instance.
\section{Conclusion and future work}\label{sec9}
In this paper, we studied the problem of integrating train timetable, rolling stock assignment and short-turning strategy for off-peak hours and with the addition of skip-stop strategy during peak hours considering a bidirectional metro line. The capacity of train, predefined operation zones, turnaround
operations, and availability of trains are all considered in the problem.
The paper presents a methodology for determining the optimal locations of depots and turnaround facilities at intermediate stations in a metro line. The calculations of various parameters such as potential train services, running time, and dwelling time are also presented in the paper. We formulated three optimization
models (each for off-peak and peak hours) with different objectives: the first model, a MILP problem, aimed at minimizing operational costs without considering passenger waiting time. Subsequently, the second model, a MINLP problem, prioritizes maximizing service quality by minimizing passenger waiting time. Lastly, we devised a multi-objective MINLP problem to simultaneously address both objectives. To improve computational efficiency, we later transformed the MINLP problem into an MILP using linearization techniques. We solved the multi-objective optimization problem using the epsilon constraint method, obtaining non-dominating solutions for both objectives. To validate our models, we tested all
the proposed formulations on simplified Santiago metro line 1. The obtained results shows that proper execution of advanced operational strategies (i.e., short-turning, rolling stock assignment and skip-stopping) improves the movement of rolling stock and offers additional train services with existing number of trains.

In this paper, we focused on medium-scale instances and time horizons. In the future, we intend to develop efficient methods to tackle larger-scale problems. Another direction for future research would be to address uncertainties related to passenger demand using a robust/stochastic optimization approach. 
\section*{Acknowledgments} We thank INFORMS Railway Applications Section for providing us with the data for the paper. An earlier version of the paper was awarded second prize in the 2023 INFORMS RAS problem solving competition.
\section*{Appendix A and Appendix B :  Supplementary data}



%

\bibliographystyle{agsm}
\newpage
\bibliography{references}

@article{barrena2014exact,
  title={Exact formulations and algorithm for the train timetabling problem with dynamic demand},
  author={Barrena, Eva and Canca, David and Coelho, Leandro C and Laporte, Gilbert},
  journal={Computers \& Operations Research},
  volume={44},
  pages={66--74},
  year={2014},
  publisher={Elsevier}
}

@article{blanco2020optimization,
  title={An optimization model for line planning and timetabling in automated urban metro subway networks. A case study},
  author={Blanco, V{\'\i}ctor and Conde, Eduardo and Hinojosa, Yolanda and Puerto, Justo},
  journal={Omega},
  volume={92},
  pages={102165},
  year={2020},
  publisher={Elsevier}
}

@article{cadarso2012integration,
  title={Integration of timetable planning and rolling stock in rapid transit networks},
  author={Cadarso, Luis and Mar{\'\i}n, {\'A}ngel},
  journal={Annals of operations research},
  volume={199},
  pages={113--135},
  year={2012},
  publisher={Springer}
}

@article{canca2018integrated,
  title={The integrated rolling stock circulation and depot location problem in railway rapid transit systems},
  author={Canca, David and Barrena, Eva},
  journal={Transportation Research Part E: Logistics and Transportation Review},
  volume={109},
  pages={115--138},
  year={2018},
  publisher={Elsevier}
}

@article{canca2014design,
  title={Design and analysis of demand-adapted railway timetables},
  author={Canca, David and Barrena, Eva and Algaba, Encarnaci{\'o}n and Zarzo, Alejandro},
  journal={Journal of Advanced Transportation},
  volume={48},
  number={2},
  pages={119--137},
  year={2014},
  publisher={Wiley Online Library}
}

@article{canca2016short,
  title={A short-turning policy for the management of demand disruptions in rapid transit systems},
  author={Canca, David and Barrena, Eva and Laporte, Gilbert and Ortega, Francisco A},
  journal={Annals of Operations Research},
  volume={246},
  pages={145--166},
  year={2016},
  publisher={Springer}
}

@article{canca2017design,
  title={Design of energy-efficient timetables in two-way railway rapid transit lines},
  author={Canca, David and Zarzo, Alejandro},
  journal={Transportation Research Part B: Methodological},
  volume={102},
  pages={142--161},
  year={2017},
  publisher={Elsevier}
}

@article{cepeda2006frequency,
  title={A frequency-based assignment model for congested transit networks with strict capacity constraints: characterization and computation of equilibria},
  author={Cepeda, Manuel and Cominetti, Roberto and Florian, Michael},
  journal={Transportation research part B: Methodological},
  volume={40},
  number={6},
  pages={437--459},
  year={2006},
  publisher={Elsevier}
}

@article{delle1998service,
  title={Service optimization for bus corridors with short-turn strategies and variable vehicle size},
  author={Delle Site, Paolo and Filippi, Francesco},
  journal={Transportation Research Part A: Policy and Practice},
  volume={32},
  number={1},
  pages={19--38},
  year={1998},
  publisher={Elsevier}
}

@article{dong2020integrated,
  title={Integrated optimization of train stop planning and timetabling for commuter railways with an extended adaptive large neighborhood search metaheuristic approach},
  author={Dong, Xinlei and Li, Dewei and Yin, Yonghao and Ding, Shishun and Cao, Zhichao},
  journal={Transportation Research Part C: Emerging Technologies},
  volume={117},
  pages={102681},
  year={2020},
  publisher={Elsevier}
}

@article{ghaemi2018macroscopic,
  title={Macroscopic multiple-station short-turning model in case of complete railway blockages},
  author={Ghaemi, Nadjla and Cats, Oded and Goverde, Rob MP},
  journal={Transportation Research Part C: Emerging Technologies},
  volume={89},
  pages={113--132},
  year={2018},
  publisher={Elsevier}
}

@article{gkiotsalitis2019cost,
  title={A cost-minimization model for bus fleet allocation featuring the tactical generation of short-turning and interlining options},
  author={Gkiotsalitis, Konstantinos and Wu, Zongxiang and Cats, Oded},
  journal={Transportation Research Part C: Emerging Technologies},
  volume={98},
  pages={14--36},
  year={2019},
  publisher={Elsevier}
}

@article{hassannayebi2016train,
  title={Train timetabling for an urban rail transit line using a Lagrangian relaxation approach},
  author={Hassannayebi, Erfan and Zegordi, Seyed Hessameddin and Yaghini, Masoud},
  journal={Applied Mathematical Modelling},
  volume={40},
  number={23-24},
  pages={9892--9913},
  year={2016},
  publisher={Elsevier}
}

@article{li2019demand,
  title={Demand-oriented train services optimization for a congested urban rail line: integrating short turning and heterogeneous headways},
  author={Li, Sijie and Xu, Ruihua and Han, Ke},
  journal={Transportmetrica A: Transport Science},
  volume={15},
  number={2},
  pages={1459--1486},
  year={2019},
  publisher={Taylor \& Francis}
}

@article{liu2020collaborative,
  title={Collaborative optimization for metro train scheduling and train connections combined with passenger flow control strategy},
  author={Liu, Renming and Li, Shukai and Yang, Lixing},
  journal={Omega},
  volume={90},
  pages={101990},
  year={2020},
  publisher={Elsevier}
}

@article{mo2021exact,
  title={An exact method for the integrated optimization of subway lines operation strategies with asymmetric passenger demand and operating costs},
  author={Mo, Pengli and D’Ariano, Andrea and Yang, Lixing and Veelenturf, Lucas P and Gao, Ziyou},
  journal={Transportation Research Part B: Methodological},
  volume={149},
  pages={283--321},
  year={2021},
  publisher={Elsevier}
}

@article{niu2013optimizing,
  title={Optimizing urban rail timetable under time-dependent demand and oversaturated conditions},
  author={Niu, Huimin and Zhou, Xuesong},
  journal={Transportation Research Part C: Emerging Technologies},
  volume={36},
  pages={212--230},
  year={2013},
  publisher={Elsevier}
}

@article{niu2015train,
  title={Train scheduling for minimizing passenger waiting time with time-dependent demand and skip-stop patterns: Nonlinear integer programming models with linear constraints},
  author={Niu, Huimin and Zhou, Xuesong and Gao, Ruhu},
  journal={Transportation Research Part B: Methodological},
  volume={76},
  pages={117--135},
  year={2015},
  publisher={Elsevier}
}

@article{tirachini2011optimal,
  title={Optimal design and benefits of a short turning strategy for a bus corridor},
  author={Tirachini, Alejandro and Cort{\'e}s, Cristi{\'a}n E and Jara-D{\'\i}az, Sergio R},
  journal={Transportation},
  volume={38},
  pages={169--189},
  year={2011},
  publisher={Springer}
}

@article{wang2018passenger,
  title={Passenger demand oriented train scheduling and rolling stock circulation planning for an urban rail transit line},
  author={Wang, Yihui and D’Ariano, Andrea and Yin, Jiateng and Meng, Lingyun and Tang, Tao and Ning, Bin},
  journal={Transportation Research Part B: Methodological},
  volume={118},
  pages={193--227},
  year={2018},
  publisher={Elsevier}
}

@article{wang2021real,
  title={Real-time integrated train rescheduling and rolling stock circulation planning for a metro line under disruptions},
  author={Wang, Yihui and Zhao, Kangqi and D’Ariano, Andrea and Niu, Ru and Li, Shukai and Luan, Xiaojie},
  journal={Transportation Research Part B: Methodological},
  volume={152},
  pages={87--117},
  year={2021},
  publisher={Elsevier}
}

@article{yalccinkaya2009modelling,
  title={Modelling and optimization of average travel time for a metro line by simulation and response surface methodology},
  author={Yal{\c{c}}{\i}nkaya, {\"O}zg{\"u}r and Bayhan, G Mirac},
  journal={European Journal of Operational Research},
  volume={196},
  number={1},
  pages={225--233},
  year={2009},
  publisher={Elsevier}
}

@article{yang2015energy,
  title={An energy-efficient scheduling approach to improve the utilization of regenerative energy for metro systems},
  author={Yang, Xin and Chen, Anthony and Li, Xiang and Ning, Bin and Tang, Tao},
  journal={Transportation Research Part C: Emerging Technologies},
  volume={57},
  pages={13--29},
  year={2015},
  publisher={Elsevier}
}

@article{yin2021timetable,
  title={Timetable coordination in a rail transit network with time-dependent passenger demand},
  author={Yin, Jiateng and D’Ariano, Andrea and Wang, Yihui and Yang, Lixing and Tang, Tao},
  journal={European Journal of Operational Research},
  volume={295},
  number={1},
  pages={183--202},
  year={2021},
  publisher={Elsevier}
}

@article{yin2017dynamic,
  title={Dynamic passenger demand oriented metro train scheduling with energy-efficiency and waiting time minimization: Mixed-integer linear programming approaches},
  author={Yin, Jiateng and Yang, Lixing and Tang, Tao and Gao, Ziyou and Ran, Bin},
  journal={Transportation Research Part B: Methodological},
  volume={97},
  pages={182--213},
  year={2017},
  publisher={Elsevier}
}

@article{yue2017integrated,
  title={Integrated train timetabling and rolling stock scheduling model based on time-dependent demand for urban rail transit},
  author={Yue, Yixiang and Han, Juntao and Wang, Shifeng and Liu, Xiang},
  journal={Computer-Aided Civil and Infrastructure Engineering},
  volume={32},
  number={10},
  pages={856--873},
  year={2017},
  publisher={Wiley Online Library}
}

@article{zhang2018short,
  title={A short turning strategy for train scheduling optimization in an urban rail transit line: the case of Beijing subway line 4},
  author={Zhang, Miao and Wang, Yihui and Su, Shuai and Tang, Tao and Ning, Bin},
  journal={Journal of Advanced Transportation},
  volume={2018},
  year={2018},
  publisher={Hindawi}
}

@article{assis2004generation,
  title={Generation of optimal schedules for metro lines using model predictive control},
  author={Assis, Wanderson O and Milani, Bas{\i}́lio EA},
  journal={Automatica},
  volume={40},
  number={8},
  pages={1397--1404},
  year={2004},
  publisher={Elsevier}
}

@article{cury1980methodology,
  title={A methodology for generation of optimal schedules for an underground railway system},
  author={Cury, J and Gomide, F and Mendes, M},
  journal={IEEE Transactions on automatic control},
  volume={25},
  number={2},
  pages={217--222},
  year={1980},
  publisher={IEEE}
}

@article{yuan2022integrated,
  title={Integrated optimization of train timetable, rolling stock assignment and short-turning strategy for a metro line},
  author={Yuan, Jiawei and Gao, Yuan and Li, Shukai and Liu, Pei and Yang, Lixing},
  journal={European Journal of Operational Research},
  volume={301},
  number={3},
  pages={855--874},
  year={2022},
  publisher={Elsevier}
}

@article{freyss2013continuous,
  title={Continuous approximation for skip-stop operation in rail transit},
  author={Freyss, Maxime and Giesen, Ricardo and Mu{\~n}oz, Juan Carlos},
  journal={Procedia-Social and Behavioral Sciences},
  volume={80},
  pages={186--210},
  year={2013},
  publisher={Elsevier}
}

@article{wang2014efficient,
  title={Efficient bilevel approach for urban rail transit operation with stop-skipping},
  author={Wang, Yihui and De Schutter, Bart and van den Boom, Ton JJ and Ning, Bin and Tang, Tao},
  journal={IEEE Transactions on Intelligent Transportation Systems},
  volume={15},
  number={6},
  pages={2658--2670},
  year={2014},
  publisher={IEEE}
}

@article{jiang2019q,
  title={Q-learning approach to coordinated optimization of passenger inflow control with train skip-stopping on a urban rail transit line},
  author={Jiang, Zhibin and Gu, Jinjing and Fan, Wei and Liu, Wei and Zhu, Bingqin},
  journal={Computers \& Industrial Engineering},
  volume={127},
  pages={1131--1142},
  year={2019},
  publisher={Elsevier}
}

@article{zhu2020integrated,
  title={Integrated timetable rescheduling and passenger reassignment during railway disruptions},
  author={Zhu, Yongqiu and Goverde, Rob MP},
  journal={Transportation Research Part B: Methodological},
  volume={140},
  pages={282--314},
  year={2020},
  publisher={Elsevier}
}

@article{chen2022real,
  title={Real-time optimization for train regulation and stop-skipping adjustment strategy of urban rail transit lines},
  author={Chen, Zebin and Li, Shukai and D’Ariano, Andrea and Yang, Lixing},
  journal={Omega},
  volume={110},
  pages={102631},
  year={2022},
  publisher={Elsevier}
}

@misc{RASCompetition2023,
  title = {RAS Problem Solving Competition},
  howpublished = {\url{https://connect.informs.org/railway-applications/new-item3/problem-solving-competition681}},
  note = {Accessed: [12.02.2024]},
  year = {2023}
}

@article{sun2014demand,
  title={Demand-driven timetable design for metro services},
  author={Sun, Lijun and Jin, Jian Gang and Lee, Der-Horng and Axhausen, Kay W and Erath, Alexander},
  journal={Transportation Research Part C: Emerging Technologies},
  volume={46},
  pages={284--299},
  year={2014},
  publisher={Elsevier}
}

@article{haimes1971bicriterion,
  title={On a bicriterion formulation of the problems of integrated system identification and system optimization},
  author={Haimes, Yacov},
  journal={IEEE transactions on systems, man, and cybernetics},
  number={3},
  pages={296--297},
  year={1971},
  publisher={Institute of Electrical and Electronics Engineers (IEEE)}
}

@book{chankong2008multiobjective,
  title={Multiobjective decision making: theory and methodology},
  author={Chankong, Vira and Haimes, Yacov Y},
  year={1983},
  publisher={Courier Dover Publications}
}






\end{document}


\noindent\textbf{\textit{Appendix A}}\\
As previously stated, certain constraints and the objective function of the model are non-linear. To address this, we must linearize these constraints and the objective function using specific following transformation.\\
\textbf{Lemma 1.} 
For $k\in K^{(up)}$, $i$, $j$ $\in$ $S^{(up)}$, where $i < j$, constraint (50) can be reformulated into the following set of linear constraints.
\begin{align*}
& 0 \leq w^{b}_{k,i,j} \leq w_{k,i,j} & (75)\\
&  w^{b}_{k,i,j} \leq M.x_{k,i} & (76)\\
&  w^{b}_{k,i,j} \leq M.x_{k,j} & (77)\\
&  w^{b}_{k,i,j} \geq w_{k,i,j} - M.(2-x_{k,i} - x_{k,j}) & (78)
\end{align*}  
Similarly, for  $k\in K^{(dn)}$, $i$, $j$ $\in$ $S^{(dn)}$, with $i < j$, constraint (59) can be linearized. We can linearized constraint (52) by introducing additional variables.\\
\textbf{Lemma 2.} For $k\in K^{(up)}$, $i \in$ $S^{(up)}$, the constraint (52) is equivalent to the following set of linear constraints.
 \begin{align*}
&  n^{b}_{k,i} \leq w^{b}_{k,i} & (79)\\
&  n^{b}_{k,i} \leq C - n_{k,i-1} + n^{a}_{k,i} & (80)\\
&  n^{b}_{k,i} \geq w^{b}_{k,i} - M.(1-\psi_{k,i}) & (81)\\
&  n^{b}_{k,i} \geq C-n_{k,i-1}+n^{a}_{k,i}-M.\psi_{k,i} & (82)
\end{align*}
Where the binary variable $\psi_{k,i}$ establishes the relationship between $w^{b}_{k,i}$ and $C-n_{k,i-1} + n^{a}_{k,i}$, i.e.,  
\begin{align*}
&\psi_{k.i} = 1, w^{b}_{k,i} \leq C-n_{k,i-1} + n^{a}_{k,i}&\label{} \\
&\psi_{k.i} = 0, w^{b}_{k,i} > C-n_{k,i-1} + n^{a}_{k,i}&\label{} 
\end{align*}
The binary variable $\psi_{k,i}$ can be accurately defined using the following linear constraints.
 \begin{align*}
 &\psi_{k.i} > \frac{(C-n_{k,i-1} + n^{a}_{k,i})-w^{b}_{k,i}}{M}& (83)\\
&\psi_{k.i} \leq 1 + \frac{(C-n_{k,i-1} + n^{a}_{k,i})-w^{b}_{k,i}}{M}&(84)\\   &\psi_{k.i} \in \{0,1\} &(85)   
\end{align*}
Similarly, for $k\in K^{(dn)}$, $i$ $\in$ $S^{(dn)}$, constraint (61) can be linearized. By introducing a continuous variable $zz_{k,m,n}$ to represent the bilinear term $(d_{k,n}-d_{k-1,m}) \cdot (z_{k,m,n})$ in objective function (72), the objective function can be linearized according to the following lemma.\\
\textbf{Lemma 3.} The objective function (72) can be effectively substituted with the following expression and additional constraints.
\begin{multline*}
 F^{*}_{obj,2} = \min~\sum_{k \in K^{(up)}}\sum_{m \in S^{(up)}/\{N-2,N\}} \sum_{n \in S^{(up)}/\{1,3\}} (zz_{k,m,n}) \\ +\sum_{k \in K^{(dn)}} \sum_{m \in S^{(dn)}/\{2N-2,2N\}} \sum_{n \in S^{(dn)}/ \{N+1,N+3\}}  (zz_{k,m,n})\\
+ \sum_{k \in K^{(up)}}\sum_{i \in S^{(up)}} (d_{k,i} - d_{k-1,i})+ 
 \sum_{k \in K^{(dn)}}\sum_{i \in S^{(dn)}} (d_{k,i} - d_{k-1,i}) 
\end{multline*}
\begin{align*}
&zz_{k,m,n} \geq (d_{k,n} - d_{k,m}) - M.(1-z_{k,m,n}) \nonumber \\ &\forall k \in K^{(up)}, m \in S^{(up)}/\{N-2,N\}, n \in S^{(up)}/\{1,3\}& (86)\\
&zz_{k,m,n} \geq (d_{k,n} - d_{k,m}) - M.(1-z_{k,m,n}) \nonumber \\& \forall k \in K^{(dn)}, m \in S^{(dn)}/\{2N-2,2N\}, n \in S^{(dn)}/\{N+1,N+3\} &(87) \\
&zz_{k,m,n} \leq M.(z_{k,m,n}) \nonumber \\ &\forall k \in K^{(up)}, m \in S^{(up)}/\{N-2,N\}, n \in S^{(up)}/\{1,3\} &(88)\\
&zz_{k,m,n} \leq M.(z_{k,m,n}) \nonumber \\& \forall k \in K^{(dn)}, m \in S^{(dn)}/\{2N-2,2N\}, n \in S^{(dn)}/\{N+1,N+3\} &(89) 
\end{align*}
\textbf{\textit{Appendix B}}
\begin{table}[htbp]
\caption*{Table 12 Pre-determined segment running time of simplified Santiago metro line 1}
\resizebox{\columnwidth}{!}{%
\begin{tabular}{llll}
\hline
\textbf{Upstream Direction} & \textbf{Running time (sec)} & \textbf{Downstream Direction} & \textbf{Running time (sec)} \\ \hline
\textit{San Pablo - Neptuno} & 44.8380 & \textit{Estación Central - Universidad de Santiago} & 46.5032 \\
\textit{Neptuno - Pajaritos} & 63.5149 & \textit{Universidad de Santiago - San Alberto Hurtado} & 40.7426 \\
\textit{Pajaritos - Las Rejas} & 50.0135 & \textit{San Alberto Hurtado - Ecuador} & 46.6832 \\
\textit{Las Rejas - Ecuador} & 46.0081 & \textit{Ecuador - Las Rejas} & 46.0081 \\
\textit{Ecuador - San Alberto Hurtado} & 46.6832 & \textit{Las Rejas - Pajaritos} & 50.0135 \\
\textit{San Alberto Hurtado - Universidad de Santiago} & 40.7426 & \textit{Pajaritos - Neptuno} & 63.5149 \\
\textit{Universidad de Santiago - Estación Central} & 46.5032 & \textit{Neptuno - San Pablo} & 44.8380 \\ \hline
\end{tabular}
}
\end{table}

\begin{table}[hbt!]
\caption*{Table 13 Pre-determined station dwelling time of simplified Santiago metro line 1}
\resizebox{\columnwidth}{!}{%
\begin{tabular}{lcclcc}
\hline
\multirow{2}{*}{\textbf{Station}} & \multicolumn{2}{c}{\textbf{Dwelling time (sec)}} & \multirow{2}{*}{\textbf{Station}} & \multicolumn{2}{c}{\textbf{Dwelling time (sec)}} \\ \cline{2-3} \cline{5-6} 
 & \multicolumn{1}{l}{\textbf{Upstream}} & \multicolumn{1}{l}{\textbf{Downstream}} &  & \multicolumn{1}{l}{\textbf{Upstream}} & \multicolumn{1}{l}{\textbf{Downstream}} \\ \hline
\textit{San Pablo} & 45 & 45 & \textit{Ecuador} & 40 & 40 \\
\textit{Neptuno} & 35 & 35 & \textit{San Alberto Hurtado} & 40 & 40 \\
\textit{Pajaritos} & 35 & 35 & \textit{Universidad de Santiago} & 35 & 35 \\
\textit{Las Rejas} & 45 & 45 & \textit{Estación Central} & 45 & 45 \\ \hline
\end{tabular}
}
\end{table}